\documentclass[12pt, leqno]{article}

\usepackage{amsmath,amsthm,amssymb,enumerate}
\usepackage{epsfig}
\usepackage{amssymb}
\usepackage{amsmath}
\usepackage{amssymb}
\usepackage{amsmath,amsthm}
\usepackage[latin1]{inputenc}
\usepackage[T1]{fontenc}
\usepackage{path}
\usepackage{ae,aecompl}
\usepackage{amsfonts}
\usepackage{amsxtra}
\usepackage{bbm,euscript,mathrsfs}
\usepackage{color}

\allowdisplaybreaks

\usepackage{amsfonts}
\usepackage{amsxtra}

\usepackage[frak,hyperref]{paper_diening}

\usepackage{enumitem}
\setenumerate{label={\rm (\alph{*})}}

\makeatletter
\long\def\unmarkedfootnote#1{{\long\def\@makefntext##1{##1}\footnotetext{#1}}}
\makeatother

\newcommand{\Div}{\divergence}
\newcommand{\ep}{\bfvarepsilon}
\newcommand{\R}{\mathbb R}
\newcommand{\N}{\mathbb N}

\newcommand{\rn}{\mathbb R^n}
\newcommand{\rmat}{\mathbb R^{n\times n}}

\DeclareMathOperator{\Ker}{Ker}

\newtheorem{defs}{Definition}[section]
\newtheorem{theorem}[defs]{Theorem}

\newtheorem{example}[defs]{Example}

\newtheorem{lemma}[defs]{Lemma}

\newtheorem{corollary}[defs]{Corollary}

\newtheorem{remark}[defs]{Remark}
\newtheorem{props}[defs]{Proposition}

\title{Negative  Orlicz-Sobolev norms  and strongly nonlinear  systems in fluid mechanics
}

\author{
  Dominic Breit\\
 {\it Mathematical Institute, LMU M\"unich}\\ {\it Theresianstra\ss e 39, 80333 M\"unich, Germany}
  \bigskip
  \\
 Andrea Cianchi \\
 {\it Dipartimento di Matematica e Informatica \lq\lq U. Dini", Universit\`a di Firenze}\\ {\it Piazza Ghiberti
27, 50122 Firenze, Italy}
}

\date{}

\numberwithin{equation}{section}

\begin{document}
\maketitle

%


%
%
%
%
%
%
\date{}
%
%
%

\begin{abstract}
We prove a version of the negative norm theorem in Orlicz-Sobolev
spaces. A study of continuity properties of the \Bogovskii -operator
between Orlicz spaces is a crucial step, of independent interest, in
our approach. Applications to the  problem of pressure
reconstruction for Non-Newtonian fluids governed by constitutive
laws, which are not necessarily of power type, are presented. A key
inequality for a numerical analysis of the underlying elliptic
system is also derived.
\end{abstract}



\unmarkedfootnote{
}

\section{Introduction}\label{intro}

Assume that $\Omega$ is a domain, namely a connected open set, in
$\R^n$, with $n \geq 2$, and let $1\leq p \leq \infty$. The negative
Sobolev norm of the distributional gradient of a function $u \in
L^1(\Omega)$ can be defined as
\begin{align}\label{negp}
\|\nabla u\|_{W^{-1,p}(\Omega, \rn)}=\sup_{\bfvarphi\in
C^\infty_0(\Omega, \R^n)} \frac{\int_\Omega u\,\Div \bfvarphi\, dx
}{\|\nabla \bfvarphi\|_{L^{p'}(\Omega , \R^{n\times n})}}\,dx\,.
\end{align}
 Here, $div$ stands for the divergence operator, and $p'= \frac p{p-1}$, the H\"older conjugate of $p$. Moreover, in \eqref{negp}, and in similar
occurrences throughout the paper, we tacitly assume that the
sumpremum is extended over all functions $\bfv$ which do not vanish
identically. Observe that the notation $\|\nabla
u\|_{W^{-1,p}(\Omega, \rn)}$ is consistent with the fact that  the
quantity  on the right-hand side of \eqref{negp} agrees with the
norm of $\nabla u$, when regarded as an element of the dual of
$W_0^{1,p'}(\Omega)$, where $W_0^{1,p'}(\Omega, \R^n)$
denotes the Sobolev space of $\R^n$-valued functions in $\Omega$
with zero traces.
\par\noindent
Definition \eqref{negp} goes back to Ne\v{c}as \cite{Ne}, who showed
that, if $\Omega$ is regular enough -- a bounded Lipschitz domain,
say -- and $1 < p < \infty$, then the $L^p(\Omega)$ norm of a
function is equivalent to the $W^{-1,p}(\Omega, \rn)$  norm of its
gradient. Namely, there exist positive constants $C_1=C_1(\Omega ,
p)$ and $C_2=C_2(n)$,
such that
\begin{align}\label{eq:1.1}
C_1\|u-u_\Omega\|_{L^p(\Omega)}\leq \|\nabla u\|_{W^{-1,p}(\Omega,
\rn)}\leq C_2\|u-u_\Omega\|_{L^p(\Omega)}
\end{align}
for every $u\in L^1(\Omega)$, where
$$u_\Omega=\frac 1{|\Omega|}\int _\Omega u\,dx$$
the mean value of $u$ over $\Omega$, and $|\Omega|$ denotes the
Lebesgue measure of $\Omega$. This result is known as Ne\v{c}as
negative norm theorem.
\par
In the present paper we are concerned with a version of this theorem
when a negative norm is introduced with the Lebesgue space
$L^p(\Omega)$
replaced with a  general  Orlicz space.
Loosely speaking,  Orlicz spaces extend  Lebesgue
spaces in that the role of the power $t^p$ in their definition is
played by a more general convex function  (a precise definition is
recalled in the next section). Clearly, a full analogue of
\eqref{eq:1.1} cannot hold for arbitrary Orlicz spaces, since
\eqref{eq:1.1} fails, for instance, in the borderline cases when
either $p=1$, or $p=\infty$. Our main result in this connection
asserts that, however, an inequality in the spirit of \eqref{eq:1.1}
still holds if, on the leftmost side, an Orlicz norm appears which,
in general, has to be slightly weaker than that on the rightmost
side. A precise balance between the relevant norms for a conclusion
of this kind to hold is the content of Theorem \ref{thm:3.2},
Section \ref{main}. A key step in our approach is an analysis, of
possible independent interest, of the divergence equation in Orlicz
spaces, via boundedness properties of the (gradient of the)
\Bogovskii \,\,operator in these spaces.

\par
Our main  motivation for a discussion of negative Orlicz norms are
applications to   a mathematical model for Non-Newtonian fluids. In
the stationary case, the relevant model tells us that the velocity
field $\bfv:\Omega\rightarrow\R^n$ and the pressure
$\pi:\Omega\rightarrow\R$ of a fluid solve the following system  of
partial differential equations:
\begin{align}\label{eq:1.1syst}
\begin{cases}
-\Div \bfS+\varrho\Div\big(\bfv\otimes\bfv\big)+\nabla \pi =\varrho\Div\bfF \quad & \hbox{in $\Omega$,}\\
\Div \bfv  =0 & \hbox{in $\Omega$,}\\
\bfv =0 &\hbox{on $\partial \Omega$,}\end{cases}
%
%
\end{align}
see, for instance, \cite{BAH}. Here, $\varrho$ is a positive
constant, whose physical meaning is the density of the fluid, the
operation $\otimes$ denotes tensor product, the function
$\bfF:\Omega\rightarrow\R^{n\times n}$ describes the given volume
forces,  and  the stress deviator $\bfS:\Omega\rightarrow\R^{n\times
n}$ is related to $\bfv$ via a constitutive law. Of course, the
physically relevant dimensions are $n=2$ and $n=3$. In a common
model for fluids with Non-Newtonian behavior, the dependence of
$\bfS$ on $\bfv$ is through a nonlinear function of the symmetric
part $\ep(\bfv)$ of its $\rmat$-valued gradient, defined as
$\ep(\bfv)=\frac{1}{2}\big(\nabla \bfv+(\nabla\bfv)^T\big)$, where
$\lq\lq (\cdot )^T "$ stands for transpose. Lebesgue and usual
Sobolev spaces provide an appropriate functional framework for the
study of existence, uniqueness and regularity of solutions to
\eqref{eq:1.1syst} when this nonlinear function is of power type. On
the other hand, if nonlinearities of non-polynomial type are
allowed, such as in the Eyring-Prandtl model \cite{E}, the more
flexible Orlicz and Orlicz-Sobolev spaces have to be called into
play. In particular, in Section \ref{fluids} we show how our
negative Orlicz norm theorem applies in the description of a
suitable space for the pressure $\pi$ in \eqref{eq:1.1syst}. A
related result of use in a finite elements method for a numerical
analysis of \eqref{eq:1.1syst} is also   established.

\section{Orlicz and Orlicz-Sobolev spaces}

In the present section we collect some definitions and fundamental
results from the theory of Orlicz and Orlicz-Sobolev spaces. We
refer to \cite{RR1, RR2} for a comprehensive treatment of this
topic.
\par
 A function $A: [0, \infty) \to [0, \infty]$ is called a
Young function if it is convex, left-continuous, and neither
identically equal to $0$, nor to $\infty$. Thus, with any such
function,  it is uniquely associated a (nontrivial) non-decreasing
left-continuous function $a:[0, \infty) \rightarrow [0, \infty]$
such that
\begin{equation}\label{B.4}
A(s) = \int _0^s a(r)\,dr \qquad {\rm for} \,\, s\geq 0.
\end{equation}
The Young conjugate $\widetilde{A}$ of $A$ is the Young function
defined by
$$\widetilde{A}(s) = \sup \{rs-A(r):\,r\geq 0\} \qquad {\rm for}\qquad  s\geq 0\,.$$
Note the representation formula
$$\widetilde A (s) = \int _0^s a^{-1}(r)\,dr \qquad {\rm for} \,\, s\geq
0,$$ where $a^{-1}$ denotes the (generalized) left-continuous
inverse of $a$. One has that
\begin{equation}\label{AAtilde}
r \leq A^{-1}(r)\widetilde A^{-1}(r) \leq 2 r \quad \hbox{for $r
\geq 0$.}
\end{equation}
Moreover,
\begin{equation}\label{B.7'}
\widetilde{\widetilde{A}}=A\,
\end{equation}
for any Young function $A$.
If $A$ is any Young function and $\lambda \geq 1$, then
\begin{equation}\label{lambdaA}
\lambda A(s) \leq A(\lambda s) \quad \hbox{for $s \geq 0$.}
\end{equation}
As a consequence, if $\lambda \geq 1$, then
\begin{equation}\label{lambdaA-1}
 A^{-1}(\lambda s) \leq \lambda A^{-1}( s) \quad \hbox{for $s \geq 0$,}
\end{equation}
where $A^{-1}$ denotes the (generalized) right-continuous inverse of
$A$.
\par\noindent
 A Young function $A$ is said to satisfy the
$\Delta_2$-condition
 if there exists a positive constant $C$ such that
\begin{align}\label{delta2}
A(2s)\leq CA(s)\quad \textrm{for \,\,} s\geq 0.
\end{align}
 We say that $A$ satisfies the $\nabla_2$-condition if
there exists a constant $C>2$ such that
\begin{equation}\label{nabla2}
A(2s) \geq C A(s)
\end{equation}
for $s \geq 0$. If  \eqref{delta2} [resp. \eqref{nabla2}] just holds
for $s \geq s_0$ for some $s_0
>0$, then $A$ is said to satisfy the  $\Delta_2$-condition [$\nabla_2$-condition] near
infinity. We shall also write $A\in \Delta _2$ [$A \in \nabla _2$]
to denote that $A$ satisfies the  $\Delta_2$-condition
[$\nabla_2$-condition].
\par\noindent
One has that $A \in \Delta_2$ [near infinity] if and only if
$\widetilde A \in \nabla_2$-condition [near infinity].
%
%
%
%
%
%
%
%
%
%
%
\par\noindent
A Young function $A$ is said to dominate another Young function $B$
[near infinity] if there exists a positive constant  $C$
\begin{equation}\label{B.5bis}
B(s)\leq A(C s) \qquad \textrm{for \,\,\,$s\geq 0$\,\, [$s\geq s_0$
\,\, for some $s_0>0$]\,.}
\end{equation}
The functions $A$ and $B$ are called equivalent [near infinity] if
they dominate each other [near infinity].
\par
Let $\Omega$ be a measurable subset of $\R^n$, and let $A$ be a
Young function. The Luxemburg norm, associated with  $A$,  is
defined as
\begin{align*}
\|u\|_{L^A(\Omega)}=\inf\left\{\lambda:\,\,\int_\Omega
A\Big(\frac{|u(x)|}{\lambda}\Big)\,dx\leq 1\right\}
\end{align*}
for any measurable function $u: \Omega \to \mathbb R$. The
collection of all  functions $u$ for which such norm is finite is
called the  Orlicz space $L^A(\Omega)$, and  is a Banach function
space. The subspace of $L^A(\Omega)$ of those functions $u$ such
that $\int _\Omega u(x)\, dx = 0$ will be denoted by $L^A_\bot
(\Omega)$.
 A H\"older type inequality in Orlicz
spaces takes the form
\begin{equation}\label{B.5}
\|v\|_{L^{\widetilde{A}}(\Omega)}\leq \sup _{u \in L^A(\Omega)}
\frac {\int _\Omega u(x) v(x)\, dx}{\|u\|_{L^A(\Omega)}} \leq 2
\|v\|_{L^{\widetilde{A}}(\Omega)}
\end{equation}
for every $v\in L^{\widetilde{A}}(\Omega)$.
 If
$A$ dominates $B$, then
\begin{equation}\label{B.6}
L^A(\Omega)\to L^B(\Omega),
\end{equation}
with embedding norm depending on the constant $C$ appearing in
\eqref{B.5bis}. When $|\Omega | < \infty$, embedding \eqref{B.6}
also holds if $A$ dominates $B$ just near infinity, but, in this
case, the embedding constant also depends on $B$, $s_0$ and
$|\Omega|$.
The decreasing rearrangement $u^{\ast}: [0, \infty ) \to [0,
\infty]$ of measurable function $u: \Omega \to \mathbb R$ is the
(unique) non-increasing, right-continuous function  which is
equimeasurable with $u$. Thus,
\begin{equation*}
u^{\ast}(s) = \sup \{t\geq 0: |\{x\in \o: |u(x)|>t \}|>s \} \qquad
\mathrm{for} \, s\geq 0.
\end{equation*}
The equimeasurability of $u$ and $u^*$ implies that
\begin{equation}\label{B.3}
\|u \|_{L^A(\Omega)} = \|u^{\ast} \|_{L^A(0,|\Omega|)}
\end{equation}
for every $u\in L^A(\Omega)$.
\par
The Lebesgue spaces $L^p(\Omega)$, corresponding to the choice
$A(t)=t^p$, if $p \in [1, \infty)$, and $A(t)= 0$ for $t \in [0, 1]$
and $A(t)= \infty$ for $t>1$, if $p=\infty$, are a basic example of
Orlicz spaces. Other customary instances of Orlicz spaces are
provided by the Zygmund spaces $L^p\log^\alpha L(\Omega)$, and by
the exponential spaces $\exp L^\beta (\Omega)$. If either $p>1$ and
$\alpha \in \R$, or $p=1$ and $ \alpha \geq 0$, then $L^p\log^\alpha
L(\Omega)$ is the Orlicz space associated with a Young function
equivalent to $t^p (\log t)^\alpha$ near infinity. Given $\beta >0$,
$\exp L^\beta (\Omega )$ denotes the Orlicz space built upon a Young
function equivalent to $e^{t^\beta}$ near infinity.
\par\noindent The Orlicz space $L^{A}(\Omega ,\R ^n)$
 of $\R^n$-valued measurable functions on $\Omega$ is defined
as $L^{A}(\Omega , \R^n)  = (L^{A}(\Omega))^n$, and is equipped with
the norm given by $\|\bfu \|_{L^{A}(\Omega , \R^n)}=\|\,|\bfu |\,
\|_{L^A(\Omega)}$ for $\bfu \in L^{A}(\Omega , \R^n)$. The  Orlicz
space $L^{A}(\Omega , \mathbb R^{n\times n})$ of $\mathbb R^{n\times
n}$ matrix-valued  measurable functions on $\Omega$ is defined
analogously.
\par
Assume now that $\Omega$ is an open set.  The Orlicz-Sobolev space
$W^{1,A}(\Omega)$  is the set of all weakly differentiable functions
in $L^A(\Omega)$
 whose  gradient also belongs
to $L^A(\Omega)$. It is a Banach space endowed with the norm
\begin{align*}
\|u\|_{W^{1,A}(\Omega)}=\|u\|_{L^A(\Omega)}+\|\nabla
u\|_{L^A(\Omega, \rn)}.
\end{align*}
We also define the subspace of $W^{1,A}(\Omega)$ of those functions
which vanish on $\partial \Omega$ as \begin{align*}
W^{1,A}_0(\Omega) = \{u \in W^{1,A}(\Omega) :\,\, & \hbox{the
continuation of $u$ by $0$ outside $\Omega$} \\  & \qquad \quad
\hbox{is weakly differentiable in $\R ^n$}\}.
\end{align*}
In the case when $A(t)=t^p$ for some $p \geq 1$, and $\partial
\Omega$ is regular enough, such definition of $W^{1,A}_0(\Omega)$
can be shown to reproduce the usual space $W^{1,p}_0(\Omega)$
defined as the closure in $W^{1,p}(\Omega )$ of the space $C^\infty
_0(\Omega )$ of smooth compactly supported functions in $\Omega$. In
general, the set of smooth bounded functions is dense in $L^A(\Omega
)$ only if $A$ satisfies the $\Delta _2$-condition (just near
infinity when $|\Omega|< \infty$), and hence, for arbitrary $A$, our
definition of $W^{1,A}_0(\Omega)$ yields a space which can be larger
than the closure of $C^\infty _0(\Omega )$ in $W^{1,A}_0(\Omega)$
even for smooth domains. On the other hand,
%
%
%
if
$\Omega$ is a Lipschitz domain,
 then $W^{1,A}_0(\Omega)= W^{1,A}(\Omega)\cap W^{1,1}_0(\Omega),$
where $W^{1,1}_0(\Omega)$ is defined as usual. Recall that an open
set $\Omega$ is called a Lipschitz domain if it is bounded and there
exists a neighborhood $\mathcal U$  of each point of $\partial
\Omega$ such that $\Omega \cap \mathcal U$ is  the subgraph of a
Lipschitz continuous function of $n-1$ variables. An open set
$\Omega$ is said to have the cone property if there exists a finite
cone $\Lambda$ such that each point of $\Omega$ is the vertex of a
finite cone  contained in $\Omega$ and congruent to $\Lambda$.
Clearly, any Lipschitz domain has the cone property, but the
converse is not true in general.
\par\noindent
A Poincar\'e type inequality in Orlicz-Sobolev spaces tells us that,
if $\Omega$ is a Lipschitz domain, then there exists a constant $C$,
depending on $n$ and on the Lipschitz constant of $\Omega$, such
that
\begin{equation}\label{poincare}
\|u- u_\Omega \|_{L^A(\Omega)}\leq C |\Omega |^{\frac 1n} \|\nabla
u\|_{L^A(\Omega,\R^n)}
\end{equation}
for every $u \in W^{1,A}(\Omega)$.
 Inequality  \eqref{poincare} is established in
\cite[Lemma 4.1]{CavCia} in the special case when $\Omega$ is a
ball. Its proof makes use of a rearrangement type inequality for the
norm $\|\nabla u\|_{L^A(\Omega )}$ which holds, in fact, for Sobolev
functions $u$  on any Lipschitz domain $\Omega$  \cite[Lemma 4.1 and
inequality (3.5)]{CP-arkiv}. The same  proof then applies to any
Lipschitz domain, and one can verify that the constant in the
resulting Poincar\'e inequality has the form claimed in
\eqref{poincare}.
\par
The Orlicz-Sobolev space $W^{1, A} (\Omega , \R^n)$ of $\R^n$-valued
functions is defined as $W^{1, A}(\Omega , \R^n) = \big(W^{1,
A}(\Omega )\big)^n$, and equipped with the norm $\|\bfu\|_{W^{1, A}
(\Omega ,\R^n)} =  \|\bfu \|_{L^A(\Omega ,\rn)} + \|\nabla
\bfu\|_{L^A(\Omega , \rmat)}$. The space $W^{1, A}_0(\Omega , \R^n)$
is defined accordingly.
%
%
%
%
%

\section{The negative norm theorem and the \Bogovskii \,\, operator}\label{main}

Let $A$ be a Young function, and let $\Omega$ be a bounded domain in
$\R^n$. We define the negative Orlicz-Sobolev norm associated with
$A$ of the distributional gradient of a function $u \in L^1(\Omega)$
as
\begin{align}\label{negorlicz}
\|\nabla u\|_{W^{-1,A}(\Omega, \rn)} =\sup_{\bfvarphi\in
C^\infty_0(\Omega, \R^n)} \frac{\int_\Omega u\,{\rm div}\, \bfvarphi
\,dx}{\|\nabla \bfvarphi\|_{L^{\widetilde{A}}(\Omega, \R^{n\times
n})}}.
\end{align}
The alternative notation $W^{-1}L^A(\Omega, \rn)$ will also
occasionally be employed to denote the negative Orlicz-Sobolev norm
$W^{-1,A}(\Omega, \rn)$ associated with the Orlicz space
$L^A(\Omega)$.
\par
Our Orlicz-Sobolev space version of the negative norm theorem
involves pairs of Young functions
  $A$ and $B$ which obey the following balance conditions:
\begin{equation}\label{1.1}
 t\int _{0} ^t \frac{ B(s)}{s^2} \,ds \leq A(ct) \qquad
\hbox{for $t \geq 0$,}
\end{equation}
and
\begin{equation}\label{1.2}
t \int _{0}^t \frac{ \widetilde A(s)}{s^2} \,ds \leq \widetilde
B(ct) \qquad \hbox{for $t \geq 0$,}
\end{equation}
for some positive constant $c$.
\par\noindent
Let us mention that assumptions \eqref{1.1} and \eqref{1.2} also
come  into play in a version of the Korn inequality for the
symmetric gradient in Orlicz spaces \cite{Ckorn}.

\begin{theorem}
\label{thm:3.2} Let $A$ and $B$ be Young functions fulfilling
\eqref{1.1} and \eqref{1.2}. Assume that $\Omega$ is a bounded
domain with the cone property in $\R ^n$, $n \geq 2$. Then there
exist  constants $C_1=C_1(\Omega , c)$ and $C_2=C_2(n)$ such that
\begin{align}\label{negineq}
C_1\|u-u_\Omega\|_{L^B(\Omega)}\leq \|\nabla
u\|_{W^{-1,A}(\Omega,\rn)}\leq C_2\|u-u_\Omega\|_{L^A(\Omega)}
\end{align}
for every $u \in L^1(\Omega)$. Here, $c$ denotes the constant
appearing in \eqref{1.1} and \eqref{1.2}.
\end{theorem}

\begin{remark}\label{3.2bis} {\rm Inequality \eqref{negineq} continues to
hold even if conditions \eqref{1.1} and \eqref{1.2} are just
fulfilled for $t \geq t_0$ for some $t_0>0$, but with constants
$C_1$ and $C_2$ depending also on $A$, $B$, $t_0$ and $|\Omega|$.
Indeed, the Young functions $A$ and $B$ can be replaced, if
necessary, with Young functions equivalent near infinity, and
fulfilling \eqref{1.1} and \eqref{1.2} for every $t>0$. Owing
\eqref{B.6}, such replacement leaves the quantities $\|\cdot
\|_{L^A(\Omega )}$,  $\|\cdot \|_{L^B(\Omega )}$ and $\|\nabla
\cdot\|_{W^{-1,A}(\Omega, \rn)}$ unchanged, up to multiplicative
constants depending on $A$, $B$, $t_0$ and $|\Omega|$.}
\end{remark}

 As recalled in Section \ref{intro}, the standard negative-norm
 Theorem expressed by \eqref{eq:1.1} breaks down in the borderline
 cases when $p=1$ or $p=\infty$.  This shows
that, in general, equation \eqref{negineq} cannot hold with $B=A$
on the left-hand side. In fact,
condition \eqref{1.1}, or
 \eqref{1.2} fails, with $B=A$, if, loosely speaking,
 the norm $\| \cdot
\|_{L^{A}(\Omega )}$ is \lq\lq close" to $\| \cdot \|_{L^{1}(\Omega
)}$, or to $\| \cdot \|_{L^{\infty}(\Omega )}$, respectively.
\par\noindent Note that, if either \eqref{1.1} or \eqref{1.2} holds, then $A$
dominates $B$ globally \cite[Proposition 3.5]{Ckorn}. In a sense,
assumptions \eqref{1.1} and \eqref{1.2} provide us with a
quantitative information on if, and how much,  the norm $\| \cdot
\|_{L^{B}(\Omega )}$ has to be weaker than  $\| \cdot
\|_{L^{A}(\Omega )}$  for a version of the negative-norm Theorem to
be restored in  Orlicz-Sobokev spaces.
\par\noindent
However, if $A\in \Delta _2$, then \eqref{1.2} certainly holds with
$B=A$  \cite[Theorem 1.2.1]{KokK}. Hence, in this case, assumption
\eqref{1.2} can be dropped in Theorem \ref{thm:3.2}.
 On the other hand, if $
A\in \nabla _2$, then $\widetilde{A}\in\Delta_2$ and hence
\eqref{1.1} holds with $B=A$, and assumption \eqref{1.1} can be
dropped in Theorem \ref{thm:3.2}. In particular, if   $A\in \Delta
_2 \cap\nabla _2$, then both conditions \eqref{1.1} and \eqref{1.2}
are fulfilled with $B=A$. Hence, we have the following corollary
which also follows from the results of \cite{DRS}.

\begin{corollary}
\label{cor:3.2}  Assume that $\Omega$ is a bounded domain with the
cone property in $\R ^n$, $n \geq 2$. Let $A$  be a Young function
in $\Delta _2 \cap \nabla _2$. Then there exist  a constants
$C=C(\Omega , A)$ and $C_2=C_2(n)$ such that
\begin{align}\label{negineqbis}
C_1\|u-u_\Omega\|_{L^A(\Omega)}\leq \|\nabla
u\|_{W^{-1,A}(\Omega,\rn)}\leq C_2\|u-u_\Omega\|_{L^A(\Omega)}
\end{align}
for every $u \in L^1(\Omega)$.
\end{corollary}

\begin{example}\label{exlog}
{\rm Assume that $A(t)$ is a Young function equivalent to $t^p \log
^\alpha (1+t)$ near infinity, where either $p>1$ and $\alpha \in\R$,
or $p=1$ and $\alpha \geq 1$. Hence, if $|\Omega|< \infty$, then
$$L^A(\Omega ) = L^p\log^\alpha L(\Omega).$$
\par\noindent
Assume that $\Omega$ is a bounded domain with the cone property in
$\R ^n$. If $p>1$, then $A\in \Delta _2\cap \nabla _2$, and hence
Corollary \ref{cor:3.2} tells us that
\begin{align}\label{exlog1}
C_1\|u-u_\Omega\|_{L^p\log^\alpha L(\Omega)}\leq \|\nabla
u\|_{W^{-1}L^p\log^\alpha L(\Omega, \rn)}\leq
C_2\|u-u_\Omega\|_{L^p\log^\alpha L(\Omega)}
\end{align}
for every $u \in L^1(\Omega)$. However, if $p=1$, then $A \in \Delta
_2$, but $A \notin \nabla _2$. An application of Theorem
\ref{thm:3.2} now yields
\begin{align}\label{exlog2}
C_1\|u-u_\Omega\|_{L \log^{\alpha -1} L(\Omega)}\leq \|\nabla
u\|_{W^{-1}L \log^\alpha L(\Omega , \rn)}\leq C_2\|u-u_\Omega\|_{L
\log^\alpha L(\Omega)}
\end{align}
for every $u \in L^1(\Omega)$.  In particular,
\begin{align}\label{exlog2}
C_1\|u-u_\Omega\|_{L ^1(\Omega)}\leq \|\nabla u\|_{W^{-1}L \log
L(\Omega , \rn)}\leq C_2\|u-u_\Omega\|_{L \log L(\Omega)}
\end{align}
for every $u \in L^1(\Omega)$. }
\end{example}

\begin{example}\label{exexp}
{\rm Let $\beta >0$, and let $A(t)$ be a Young function equivalent
to $\exp(t^\beta)$ near infinity. Then
$$L^A(\Omega ) = \exp L^\beta (\Omega)$$
if $|\Omega|< \infty$. One has  that $A\in \nabla _2$, but $A \notin
\Delta _2$. Theorem \ref{thm:3.2} ensures that, if $\Omega$ is a
bounded domain with the cone property in $\R ^n$, then
\begin{align}\label{exlog3}
C_1\|u-u_\Omega\|_{\exp  L^{\frac{\beta}{\beta +1}}(\Omega)}\leq
\|\nabla u\|_{W^{-1} \exp L^\beta(\Omega , \rn)}\leq
C_2\|u-u_\Omega\|_{\exp L^\beta (\Omega)}
\end{align}
for every $u \in L^1(\Omega)$. Moreover,
\begin{align}\label{exlog4}
C_1\|u-u_\Omega\|_{\exp  L(\Omega)}\leq \|\nabla u\|_{W^{-1}
L^\infty(\Omega , \rn)}\leq C_2\|u-u_\Omega\|_ {L^\infty (\Omega)}
\end{align}
for every $u \in L^1(\Omega)$. }
\end{example}

\bigskip

Our proof of Theorem \ref{thm:3.2} relies upon an analysis of
 the divergence equation
\begin{align}\label{eq:3.1}
\begin{cases}\Div\bfu= f & \quad\text{in }\Omega,\\
\bfu=0 & \text{on }\,\partial \Omega,\end{cases}
\end{align}
in Orlicz spaces. This is the objective of the next result. In what
follows, we set
$$ C^\infty_{0,\bot} (\Omega) = \{u \in C^\infty _0(\Omega ):
u_\Omega =0\}$$ and
$$L^A_\bot (\Omega) =  \{u \in L^A(\Omega ):
u_\Omega =0\}.$$
%
%
%
%

\begin{theorem}
\label{thm:3.1}  Assume that $\Omega$ is a bounded domain with the
cone property in $\R ^n$, $n \geq 2$. Let $A$ and $B$ be Young
functions fulfilling \eqref{1.1} and \eqref{1.2}. Then there exists
a bounded linear operator
\begin{equation}\label{bogov}\mathcal B_\Omega : L^A_\bot (\Omega) \to
W^{1,B}_0(\Omega, \R^n)
\end{equation} such that
\begin{equation}\label{bogovsmooth}\mathcal B_\Omega : C^\infty_{0,\bot} (\Omega) \to
C^\infty_0(\Omega, \R^n)
\end{equation}
and
\begin{equation}\label{div}
\Div(\mathcal B_\Omega f)= f\quad\text{in }\,\Omega\,
\end{equation}
 for every $f \in L^A_\bot (\Omega)$.
In particular, there exists a constant $C=C(\Omega, c)$ such that
\begin{equation}\label{grad}
\|\nabla (\mathcal B_\Omega f)\|_{L^B(\Omega , \rmat)} \leq C
\|f\|_{L^A(\Omega)}
\end{equation}
  and
\begin{equation}\label{gradint}
\int _\Omega B(|\nabla (\mathcal B_\Omega f)|)\, dx \leq  \int
_\Omega A(C|f|)\, dx
\end{equation}
for every $f \in L^A_\bot (\Omega)$. Here, $c$ denotes the constant
appearing in \eqref{1.1} and \eqref{1.2}.
%
\end{theorem}

The proof of Theorem \ref{thm:3.1} in turn makes use of a
rearrangement estimate, which extends those of \cite[Theorem
16.12]{BennettRudnick} and \cite{BagbyKurtz}, for a class of
singular integral operators of the form \begin{equation}\label{T} Tf
(x) = \lim _{\varepsilon \to 0^+} \int _{\{y:
|y-x|>\varepsilon\}}K(x, x-y) f(y)\,dy \quad \hbox{for $x \in
\R^n$,}\end{equation} for an integrable function $f : \rn \to \R$.
Here, the kernel  $K : \R ^n \times \R ^n \to \R$ fulfils the
following properties: (i) \begin{equation}\label{P1} K(x, x -
\lambda y ) = \lambda ^{-n} K(x, x-y) \quad \hbox{for $x, y \in
\R^n$;}
\end{equation}
(ii) \begin{equation}\label{P2} \int _{\mathbb S ^{n-1}} K(x, x - y)
\, d\mathcal H ^{n-1}(y) = 0 \quad \hbox{for $x\in \R^n$;}
\end{equation}
(iii) For every $\sigma \in [1, \infty )$, there exists a constant
$C_1$ such that
\begin{equation}\label{P3}
\bigg(\int _{\mathbb S ^{n-1}} |K(x, x - y)|^\sigma \, d\mathcal H
^{n-1}(y)\bigg)^{\frac 1{\sigma}} \leq C_1 (1 + |x|)^n \quad
\hbox{for $x\in \R^n$,}
\end{equation}
where $\mathbb S ^{n-1}$ denotes the unit sphere, centered at $0$,
in $\rn$, and $\mathcal H ^{n-1}$ stands for the $(n-1)$-dimensional
Hausdorff measure;
\par\noindent (iv)
There exists a constant $C_2$ such that
\begin{equation}\label{P4}
|K(x,y)| \leq C_2 \frac{(1 + |x|)^n}{|x-y|^n} \quad \hbox{for $x, y
\in \R^n$, $x \neq y$,}
\end{equation}
and, if $2|x-z| < |x-y|$, then
\begin{equation}\label{P5}
|K(x,y) - K(z,y)| \leq C_2 (1 + |y|)^n \frac{|x-z|}{|x-y|^{n+1}},
\end{equation}
\begin{equation}\label{P6}
|K(y,x) - K(y,z)| \leq C_2 (1 + |y|)^n \frac{|x-z|}{|x-y|^{n+1}}.
\end{equation}

\begin{theorem}
\label{singular} 
Let $\Omega$ be a bounded open set in $\R ^n$, and let $K(x,y)$ be
a kernel satisfying \eqref{P1}--\eqref{P6}. If $f \in L^1(\R ^n)$
and $f=0$ in $\R ^n \setminus \Omega$, then the singular integral
 operator $T$ given by \eqref{T}
is well defined for a.e. $x \in \R ^n$, and there exists a constant
$C=C(C_1, C_2, n, {\rm diam}(\Omega))$ such that
\begin{equation}\label{rearrestim}
(Tf)^*(s) \leq C \bigg(\frac 1s\int _0^s f^*(r)\, dr + \int
_s^{|\Omega|}f^*(r)\, \frac {dr}r \bigg) \quad \hbox{for $s \in (0,
|\Omega|)$.}
\end{equation}
\end{theorem}

As a consequence of Theorem \ref{singular}, the boundedness of
singular integral operators given by \eqref{T} between Orlicz spaces
associated with Young functions $A$ and $B$ fulfilling \eqref{1.1}
and \eqref{1.2} can be established.

\begin{theorem}
\label{singularcor} Let $\Omega$, $K$  and $T$ be as in
 Theorem \ref{singular}. Assume that $A$ and $B$ are Young functions
 satisfying \eqref{1.1} and \eqref{1.2}.
Then there exists a constant $C=C(C_1, C_2, n, {\rm diam}(\Omega),
c)$ such that
\begin{equation}\label{singularorlicz}
\|Tf\|_{L^B(\Omega)} \leq C \|f\|_{L^A(\Omega)},
\end{equation}
and
\begin{equation}\label{singularintegr}
\int _\Omega B(|Tf|)\,dx \leq \int _\Omega A(C|f|)\, dx
\end{equation}
for every $f \in L^A(\Omega)$. Here, $c$ denotes the constant
appearing in \eqref{1.1} and \eqref{1.2}.
\end{theorem}

\smallskip
\par\noindent
{\bf Proof}. By \cite[Lemma 1]{Ci2}, we have that, if $A$ and $B$
are Young functions satisfying \eqref{1.1}, then there exists a
constant $C=C(c)$ such that
\begin{equation}\label{hardy1}
\bigg\|\frac 1s\int _0^s \varphi(r)\, dr\bigg\|_{L^B(0, \infty )}
\leq C \|\varphi\|_{L^A(0, \infty)}
\end{equation}
for every $\varphi \in L^A(0, \infty)$. Moreover, if  $A$ and $B$
fulfill \eqref{1.2}, then there exists a constant $C=C(c)$ such that
\begin{equation}\label{hardy2}
\bigg\|\int _s^{\infty}\varphi(r)\, \frac {dr}r\bigg\|_{L^B(0,
\infty )} \leq C \|\varphi\|_{L^A(0, \infty)}
\end{equation}
for every $\varphi \in L^A(0, \infty)$. Combining
\eqref{rearrestim}, \eqref{hardy1} and \eqref{hardy2}, and making
use of property \eqref{B.3} yield inequality \eqref{singularorlicz}.
\par\noindent
As far as \eqref{singularintegr} is concerned, observe that,
inequalities \eqref{1.1} and \eqref{1.2} continue to hold, with the
same constant $c$, if $A$ and $B$ are replaced with $k A$ and $k B$,
where $k$ is any positive constant.
Thus, inequality \eqref{singularorlicz} continues to hold, with the
same constant $C$, after this replacement, whatever $K$ is, namely
\begin{equation}\label{singlambda}
\|Tf\|_{L^{k B}(\Omega)} \leq C \|f\|_{L^{k A}(\Omega)}
\end{equation}
for every $f \in L^{A}(\Omega)$.
 Now, given any such $f$, choose
$k =\frac{1}{ \int _\Omega A(|f|)\, dx}$. The very definition of
Luxemburg norm tells us that $\|f\|_{L^{k A}(\Omega)} \leq 1$.
Hence, by \eqref{singlambda}, $\|Tf\|_{L^{k B}(\Omega)} \leq C$. The
definition of Luxemburg norm again implies that $\int _\Omega k\,
B\big(\frac{|Tf|}{C}\big)\, dx \leq 1$, namely
\eqref{singularintegr}. \qed

\smallskip
\par\noindent
{\bf Proof of Theorem \ref{singular}}. Let $R>0$ be such that
$\Omega \subset B_R(0)$, the ball centered at $0$, with radius $R$.
Fix a smooth function $\eta : [0, \infty) \to [0, \infty)$ such that
$\eta = 1$ in $[0,3R]$ and $\eta = 0$ in $[4R, \infty)$. Define
$$\widehat{K}(x,y)=\eta(|x|)K(x,y) \quad \hbox{for $x, y \in \R^n$.}$$
By properties \eqref{P1}--\eqref{P6} of $K(x,y)$, one has that:
\begin{equation}\label{P1'}
\widehat{K}(x, x - \lambda y ) = \lambda ^{-n} \widehat{K}(x, x-y)
\quad \hbox{for $x, y \in \R^n$;}
\end{equation}
\begin{equation}\label{P2'}
\int _{\mathbb S ^{n-1}} \widehat{K}(x, x - y) \, d\mathcal H
^{n-1}(y) = 0 \quad \hbox{for $x\in \R^n$;}
\end{equation}
for every $\sigma \in [1, \infty )$, there exists a constant
$\widehat C_1=\widehat C_1(C_1, \sigma , R, n)$ such that
\begin{equation}\label{P3'}
\bigg(\int _{\mathbb S ^{n-1}} |\widehat{K}(x, x - y)|^\sigma \,
d\mathcal H ^{n-1}(y)\bigg)^{\frac 1{\sigma}} \leq \widehat C_1
\quad \hbox{for $x\in \R^n$,}
\end{equation}
where $C_1$ is the constant appearing in \eqref{P3}; there exists a
constant $\widehat C_2=\widehat C_2(C_2, R, n)$ such that
\begin{equation}\label{P4'}
|\widehat{K}(x,y)| \leq  \frac{\widehat C_2}{|x-y|^n} \hbox{for $x,
y \in \R^n$, $x \neq y$,}
\end{equation}
and, if $x \in \R^n$, $y \in \Omega$ and  $2|x-z| < |x-y|$, then
\begin{equation}\label{P5'}
|\widehat{K}(x,y) - \widehat{K}(z,y)| \leq  \widehat C_2
\frac{|x-z|}{|x-y|^{n+1}},
\end{equation}
\begin{equation}\label{P6'}
|\widehat{K}(y,x) - \widehat{K}(y,z)| \leq \widehat C_2
\frac{|x-z|}{|x-y|^{n+1}},
\end{equation}
where $C_2$ is the constant appearing in \eqref{P4}--\eqref{P6}.
\par\noindent
Define
\begin{align*}
\widehat{T}_{\varepsilon}f(x)&=\int_{\{y: |y-x|>\varepsilon\}}\widehat{K}_{\varepsilon}(x,y)f(y)\,dy,\\
\widehat{T}_Sf(x)& =\sup_{\varepsilon >0}|\widehat{T}_{\varepsilon}(f)(x)|.
\end{align*}
Inequality \eqref{rearrestim} will follow if we prove that
\begin{equation}\label{rearrmax}
(\widehat{T}_Sf)^*(s) \leq C \bigg(\frac 1s\int _0^s f^*(r)\, dr +
\int _s^{|\Omega|}f^*(r)\, \frac {dr}r \bigg) \quad \hbox{for $s \in
(0, \infty)$}
\end{equation}
for some constant $C=C(C_1, C_2, n, R)$, and for every $f \in L^1(\R
^n)$ such that $f =0$ in $\R ^n \setminus B_R(0)$. A proof of
inequality \eqref{rearrmax} can be accomplished along the same lines
as that of Theorem 1 of \cite{BagbyKurtz}, which in turn relies upon
similar techniques as in \cite{CoifmanFefferman}. For completeness,
we sketch such  proof hereafter.
\par\noindent
The key step in the derivation of \eqref{rearrmax} consists in
showing that, for every $\gamma \in (0,1)$, there exists a constant
$C=C(C_1, C_2, \gamma , n, R)$ such that
\begin{equation}\label{BK}
(\widehat{T}_Sf)^*(s) \leq C (M f)^*(\gamma s) +
(\widehat{T}_Sf)^*(2s) \quad \hbox{for $s \in (0, \infty)$}
\end{equation}
for every $f \in L^1(\R ^n)$ such that $f =0$ in $\R ^n \setminus
B_R(0)$. Fix $s >0$, and define
$$E= \{x \in \R ^n: \widehat{T}_Sf(x) > (\widehat{T}_Sf)^*(2s)\}.$$
Then, there exists an open set $U \supset E$ such that $|U| \leq
3s$. By Whitney's covering theorem, there exist a  family of
disjoint cubes $\{Q_k\}$ such that $U = \cup _{k=1}^\infty Q_k$,
$\sum _{k=1}^\infty |Q_k| = |U| \leq 3s$, and
$${\rm diam}(Q_k) \leq {\rm dist}(Q_k, \R^n \setminus U) \leq 4 {\rm
diam}(Q_k) \quad \hbox{for $k \in \N$.}$$ The operator
$\widehat{T}_S$ is of weak type (1,1), namely, there exists a
constant $C'$ such that
\begin{equation}\label{weak}
|\{x \in \R ^n: \widehat{T}_S f(x) > \lambda \}| \leq \frac {C'}
\lambda \|f\|_{L^1(\R ^n)}
\end{equation}
for $f \in L^1(\R ^n)$, as proved in \cite[Proof of Lemma
6.3]{DieningRuzicka}.
\par\noindent We shall now show that there
exists a constant $\overline C$
%
such that
\begin{equation}\label{3.1}
|\{x \in Q_k: \widehat{T}_S f(x) > \overline C Mf(x) +
(\widehat{T}_Sf)^*(2s)\}| \leq \frac{1-\gamma} 3 |Q_k| \quad
\hbox{for $k \in \N$.}
\end{equation}
Fix any $k \in \N$, choose $x_k \in \R ^n \setminus U$ such that
${\rm dist} (x_k , Q_k) \leq 4 {\rm diam }(Q_k),$
 and denote by $Q$ the cube, centered at $x_k$, with ${\rm diam
 }(Q)= 20 {\rm diam }(Q_k)$. Define
 $$g = f \chi _Q, \quad h = f \chi _{\R ^n \setminus Q},$$
 so that $f = g +h$. If we prove that there exist constants $\overline C_1$
 and $\overline C_2$ such that
 \begin{equation}\label{a}
\widehat{T}_S h(x) \leq \overline C_1 Mf(x) + (\widehat{T}_Sf)^*(2s)
\quad \hbox{for $x \in Q_k$,}
\end{equation}
and
\begin{equation}\label{b}
|\{x \in Q_k: \widehat{T}_S g(x) >\overline C_2 Mf(x)\}| \leq
\frac{1-\gamma} 3 |Q_k|,
\end{equation}
then \eqref{3.1} follows with $\overline C=\overline C_1 + \overline
C_2$. Consider \eqref{b} first. Let $\overline C_2$
be such that $\frac{C' |Q|}{\overline C_2} \leq \frac{1-\gamma} 3
|Q_k|$. Let $\lambda = \frac{\overline C_2}{|Q|}\int_{Q} |g|dx$.
Since $\overline C_2 Mf(x) \geq \lambda$ for $x \in Q_k$, an
application of \eqref{weak} with this choice of $\lambda$ tells us
that
\begin{align*}|\{x \in Q_k: \widehat{T}_S
g(x) >\overline C_2 Mf(x)\}| & \leq |\{\widehat{T}_S g(x) >\lambda\}| \\
& \leq \frac {C'} \lambda \int_{Q} |g|dx \leq \frac{C'
|Q|}{\overline C_2} \leq \frac{1-\gamma} 3 |Q_k|,
\end{align*}
namely \eqref{b}. In order to establish \eqref{a}, it suffices to
prove that,
 for every
$\varepsilon
>0$,
\begin{equation}\label{1}
|\widehat{T}_{\varepsilon} h(x) | \leq  \overline C_1 Mf(x) +
\widehat{T}_S f(x_k) \quad \hbox{for $x \in Q_k$.}
\end{equation}
Indeed, since $x_k \notin U$, we have that $\widehat{T}_S f(x_k)
\leq (\widehat{T}_Sf)^*(2s)$, and hence \eqref{1} implies \eqref{a}.
We may thus focus on \eqref{1}. Fix $\varepsilon >0$, and set $r=
\max \{\varepsilon , {\rm dist}(x_k , \R^n \setminus Q)\}$. Observe
that $r > 10 \,{\rm diam}(Q_k)$. Given any $x \in Q_k$, define $V =
B_\varepsilon (x) \triangle B_\varepsilon (x_k)$. One has that
\begin{align}\label{2}
|\widehat{T}_{\varepsilon} h(x) | & = \bigg|\int _{\{y:
|y-x|>\varepsilon\}} \widehat{K}(x,y) h(y)\, dy\bigg| \\ \nonumber &
\leq \bigg|\int _{\{y: |y-x_k|>\varepsilon\}} \widehat{K}(x,y)
h(y)\, dy\bigg| + \int _V | \widehat{K}(x,y) h(y)|\, dy.
\end{align}
Observe  that, if $y \in {\rm supp}\, h$, then  $|x-y| > \tfrac r2$
and hence $\frac 1{|x-y|^n} < \frac {2^n}{r^n}$. Thus, owing to
\eqref{P4'},
$$|\widehat{K}(x,y)| \leq \frac{\widehat C_2}{r^n}.$$
Moreover, $V \subset B_{3r}(x)$. Therefore, there exists a constant
$\widehat C$ such that
\begin{align}\label{3}
 \int _V | \widehat{K}(x,y) h(y)|\, dy & \leq \frac{\widehat C}{|B_{3r}(x)|}
 \\ \nonumber &\int_{B_{3r}(x)} |h(y)|\, dy \leq \widehat C M h(x) \leq \widehat C M
 f(x).
 \end{align}
 On the other hand,
 \begin{align}\label{4}
\bigg|\int _{\{y: |y-x_k|>\varepsilon\}} & \widehat{K}(x,y) h(y)\,
dy\bigg|  \leq \bigg|\int _{\{y: |y-x_k|>r\}} \widehat{K}(x,y)
h(y)\, dy\bigg|
\\ \nonumber &
\leq \bigg|\int _{\{y: |y-x_k|>r\}} \widehat{K}(x_k,y) f(y)\,
dy\bigg|
 + \int _{\{y: |y-x_k|>r\}} |\widehat{K}(x_k,y)- \widehat{K}(x,y)| \,|f(y)|\, dy,
\\ \nonumber &
\leq \widehat{T}_S (x_k)
 + \int _{\{y: |y-x_k|>r\}} |\widehat{K}(x_k,y)- \widehat{K}(x,y)| \,|f(y)|\, dy,
\end{align}
where the first inequality holds since $h(y)=0$ in $\{y: |y-x_k|\leq
r\}$ if $r= {\rm dist}(x_k , \R^n \setminus Q)$, and trivially holds
(with equality) if $r=\varepsilon$. Since $2|x-x_k| \leq |x-y|$ in
the last integral in \eqref{4}, and $f$ vanishes in $\R ^n \setminus
B_R(0)$, by \eqref{P5'}
$$|\widehat{K}(x_k,y)- \widehat{K}(x,y)| \leq \widehat C_2 \frac{|x_k - x|}{|x-y|^{n+1}}
\leq \widehat C_2  \frac{{\rm diam}(Q_k)}{|x-y|^{n+1}}.$$ Hence,
\begin{align}\label{5}
\begin{aligned}
\int _{\{y: |y-x_k|>r\}} &|\widehat{K}(x_k,y)- \widehat{K}(x,y)|
\,|f(y)|\, dy \\&\leq \int _{\{y: |y-x|>{\rm diam}(Q_k)\}}
|f(y)|\frac{{\rm diam}(Q_k)}{|x-y|^{n+1}}\, dy \leq \widetilde C
Mf(x)
\end{aligned}
\end{align}
for some constant $\widetilde C$. Note that, in the first
inequality, we have made use of the inclusion $\{y: |y-x_k|>r\}
\subset \{y: |y-x|>{\rm diam}(Q_k)\}$, which holds since $|x-x_k|<
5\, {\rm diam}(Q_k)$, and $10\, {\rm diam}(Q_k) < r$.
\par\noindent
Combining inequalities \eqref{2}-- \eqref{5} yields \eqref{1}.
Inequality \eqref{3.1} is fully established. Via summation in $k \in
Q_k$, we obtain from \eqref{3.1} that
\begin{equation}\label{3.1sum}
|\{x \in \R ^n: \widehat{T}_S f(x) > \widehat C Mf(x) +
(\widehat{T}_Sf)^*(2s)\}| \leq (1-\gamma) s\,.
\end{equation}
Coupling \eqref{3.1sum} with the inequality
\begin{equation}\label{6}
|\{x \in \R ^n:  Mf(x) > (Mf)^*(\gamma s)\}| \leq \gamma s\,
\end{equation}
tells us that
\begin{align*}
|&\{x \in \R ^n: \widehat{T}_S f(x) > \widehat C (Mf)^*(\gamma s) +
(\widehat{T}_Sf)^*(2s)\}| \leq |\{x \in \R ^n: \widehat{T}_S f(x) >
\widehat C Mf(x) + (\widehat{T}_Sf)^*(2s)\}|\\& + |\{x \in \R ^n:
Mf(x)
> (Mf)^*(\gamma s)\}|  \leq s\,,
\end{align*}
whence \eqref{BK} follows, by the very definition of decreasing
rearrangement.
\par\noindent
Starting from  inequality \eqref{BK} leads to \eqref{rearrestim},
via the same iteration argument as in the proof of \cite[Theorem
1]{BagbyKurtz}. \qed

\begin{lemma}\label{galdi}
Let $\Omega$ be a bounded domain with the cone property in $\R^n$,
with $n \geq 2$. Then there exist $N \in \N$ and  a finite family
$\{\Omega_i\}_{i=0, \dots N}$ of domains which are starshaped with
respect to balls, such that $\Omega = \cup _{i=0}^N \Omega_i$.
Moreover, given $f \in   L^A_\bot (\Omega)$, there exist $f_i \in
L^A_\bot (\Omega)$, $i=0, \dots N$, such that $f_i = 0$ in $\Omega
\setminus \Omega_i$,
%
$$ f= \sum _{i=0}^N f_i$$
and \begin{equation}\label{constant} \|f_i\|_{L^A(\Omega)}\leq C
\|f\|_{L^A(\Omega)} \quad \hbox{for $i=0, \dots ,N$,}
\end{equation}
 for some constant $C=C(\Omega)$.
\end{lemma}

\smallskip
\par\noindent
{\bf Proof, sketched}. Any bounded open set with the cone property
can be decomposed into a finite union of Lipschitz domains
\cite[Lemma 4.22]{Ad}. On the other hand, any Lipschitz domain can
be decomposed into a finite union of open sets which are starshaped
with respect to balls \cite[Lemma 3.4, Chapter 3]{Ga}. This proves
the existence of the domains $\{\Omega_i\}_{i=0, \dots N}$ as in the
statement.
 The same argument as in the
proof of \cite[Lemma 3.2, Chapter 3]{Ga} then enables one to
construct the desired   family of functions   $f_i$ on $\Omega$,
$i=1, \dots , N$, according to the following iteration scheme. We
set $G_i = \cup _{j=i+1}^N \Omega _j$, $g_0 = f$, and, for $i=1,
\dots , N-1$,
\begin{equation}\label{gk}
g_i(x) = \begin{cases} \big(1 -\chi_{\Omega _i \cap
G_i}(x)\big)g_{i-1}(x) - \frac{\chi_{\Omega _i \cap G_i}(x)}{|\Omega
_i \cap G_i|}\int _{G_i \setminus \Omega _i} g_{i-1}(y)dy & \hbox{if
$x \in G_i$,}
\\
0 & \hbox{otherwise},
\end{cases}
\end{equation}
and
\begin{equation}\label{fk}
f_i(x) = \begin{cases} g_{i-1}(x) -  \frac{\chi_{\Omega _i \cap
G_i}(x)}{|\Omega _i \cap G_i|}\int _{ \Omega _i} g_{i-1}(y)dy &
\hbox{if $x \in \Omega _i$,}
\\
0 & \hbox{otherwise}.
\end{cases}
\end{equation}
Observe that, since $\Omega$ is connected, we can always relabel the
sets $\Omega_i\cap G_i$ in such a way that $|\Omega _i \cap G_i|
>0$ for $i=1, \dots , N-1$. Finally, we
define
\begin{equation}\label{fN}
f_N = g_{N-1}.
\end{equation}
The family $\{f_i\}$ satisfies the required properties. The only
nontrivial one is \eqref{constant}. To verify the latter, fix $i$,
and observe that, by \eqref{fk}, the second inequality in
\eqref{B.5}, inequality \eqref{AAtilde}, and inequality
\eqref{lambdaA-1}
\begin{align}\label{constant1}
\|f_i\|_{L^A(\Omega )} & \leq \|g_{i-1}\|_{L^A(\Omega )}\bigg(1 +
\frac 2{|\Omega _i \cap G_i|}\|1\|_{L^A(\Omega _i \cap G_i)}
\|1\|_{L^{\widetilde A}(\Omega _i)}\bigg)
\\ \nonumber
&  = \|g_{i-1}\|_{L^A(\Omega )}\bigg(1 + \frac 2{|\Omega _i \cap
G_i|A^{-1}(1/|\Omega _i \cap G_i|)}\frac 1{\widetilde
A^{-1}(1/|\Omega _i|)}\bigg)
\\ \nonumber
 & \leq \|g_{i-1}\|_{L^A(\Omega )}\bigg(1 +
4\frac {\widetilde A^{-1}(1/|\Omega _i \cap G_i|)}{\widetilde
A^{-1}(1/|\Omega _i|)}\bigg)
\\ \nonumber
 & \leq \|g_{i-1}\|_{L^A(\Omega )}\bigg(1 +
4\frac {|\Omega _i|}{|\Omega _i \cap G_i|}\bigg).
\end{align}
On the other hand,  by \eqref{gk} and a  chain similar to
\eqref{constant1}, one has that
\begin{align}\label{constant2}
\|g_{i-1}\|_{L^A(\Omega )} & \leq \|g_{i-2}\|_{L^A(\Omega )}\bigg(1
+ \frac 2{|\Omega _{i-1} \cap G_{i-1}|}\|1\|_{L^A(\Omega _{i-1} \cap
G_{i-1})} \|1\|_{L^{\widetilde A}(G_{i-1})}\bigg)
\\ \nonumber &
\leq \|g_{i-2}\|_{L^A(\Omega )}\bigg(1 +  4\frac{\widetilde
A^{-1}(|\Omega _{i-1} \cap G_{i-1}|)}{\widetilde
A^{-1}(|G_{i-1}|)}\bigg)
\\ \nonumber & \leq \|g_{i-2}\|_{L^A(\Omega )}\bigg(1 +  4\frac{\widetilde
A^{-1}(|\Omega _{i-1} \cap G_{i-1}|)}{\widetilde A^{-1}(|G_{i-1}
|)}\bigg)
\\ \nonumber & \leq \|g_{i-2}\|_{L^A(\Omega )}\bigg(1 +  4
\max\bigg\{1, \frac{|G_{i-1}|}{|\Omega _{i-1} \cap
G_{i-1}|}\bigg\}\bigg).
\end{align}
From \eqref{constant1}, and an iteration of \eqref{constant2}, one
infers that
\begin{align}\label{constant3}
\|f_i\|_{L^A(\Omega )}  \leq \bigg(1 + 4\frac {|\Omega _i|}{|\Omega
_i \cap G_i|}\bigg) \prod_{j=1}^{i-1}\bigg(1 +  4 \max\bigg\{1,
\frac{|G_{j} |}{|\Omega _{j} \cap
G_{j}|}\bigg\}\bigg)\|f\|_{L^A(\Omega )},
\end{align}
and \eqref{constant} follows. \qed

\medskip
\par\noindent
{\bf Proof of Theorem \ref{thm:3.1}}
By Lemma \ref{galdi}, it suffices to prove the statement in the case
when $\Omega$ is a domain  starshaped with respect to a ball $B$,
which, without loss of generality, can be assumed to be centered at
the origin and with radius $1$. In this case, we shall show that the (gradient of the)
Bogovskii operator $\mathcal B_\Omega$, defined at a function $f \in
L^A_\bot (\Omega )$ as
\begin{align}\label{bog}
\mathcal B_\Omega f(x)=\int_{\Omega}
f(y)\Big(\frac{x-y}{|x-y|^{n}}\int_{|x-y|}^\infty\omega\bigg(y+\zeta\frac{x-y}{|x-y|}\bigg)\zeta^{n-1}\,d\zeta\Big)\,dy
\quad \hbox{for $x \in \Omega$,}
\end{align}
where $\omega$ is any (nonnegative) function in $C^\infty_0(B)$ with
$\int_{B} \omega\,dx=1$, agrees with a singular integral operator,
whose kernel fulfills \eqref{P1}--\eqref{P6}, plus two operators
enjoying stronger boundedness properties.  To be more precise, we
set
 $\bfu = \mathcal B_\Omega f$ and claim that $\bfu \in W^{1,1}_0(\Omega,\R^n)$, and
\begin{equation}\label{partialbog}
\frac{\partial u_i}{\partial x_j} = H_{ij} f \quad \hbox{for a.e. $x
\in \Omega$,}
 \end{equation}
 where $H_{ij}$ is the linear operator defined at
$f$ as
\begin{align}\label{H}
(H_{ij} f)
(x)&=\int_{\Omega}K_{ij}(x,y)f(y)\,dy+\int_{\Omega}G_{ij}(x,y)f(y)\,dy\\
&\quad
+f(x)\int_\Omega\frac{(x-y)_i(x-y)_j}{|x-y|^2}\omega(y)\,dy\quad
\,\, \hbox{for  $x \in \Omega$,}
\end{align}
for $i, j =1, \dots n$. Here, $K_{ij}$ is the kernel of a singular
integral operator satisfying the same assumptions as the kernel $K$
in Theorem \ref{singular}, and the kernels $G_{ij}$ satisfy
\begin{align}\label{G}
|G_{ij}(x,y)|\leq\frac{c}{|x-y|^{n-1}} \quad \hbox{for $x, y \in
\R^n$,\, $x \neq y$.}
\end{align}
To verify this assertion, recall that, if $f \in C^\infty_{0,
\bot}(\Omega)$, then $\bfu \in C^\infty _0(\Omega,\R^n)$, and
moreover equations \eqref{partialbog} and \eqref{div} hold for every
$x \in \Omega$ \cite[Proof of Lemma III.3.1]{Ga}. Consider next the
general case when $f \in L^A_\bot (\Omega )$. Owing to \eqref{1.1},
$L^A_\bot (\Omega) \to L{\rm Log} L _\bot (\Omega)$, since $B(t)$
grows at least linearly near infinity, and hence $A(t)$ dominates
the function $t \log (1+t)$ near infinity. Since the space
$C^\infty_{0, \bot}(\Omega)$ is dense in $L{\rm log} L _\bot
(\Omega)$, there exists a sequence of functions $\{f_k\} \subset
C^\infty_{0, \bot}(\Omega)$ such that $f_k \to f$ in $L{\rm log} L
(\Omega)$. One has that
$$\mathcal B_\Omega : L{\rm log} L (\Omega) \to L^1(\Omega , \rn)$$
(in fact, $\mathcal B_\Omega$ is also bounded into $L{\rm log} L
(\Omega,\R^n)$). Furthermore,
$$H_{ij} : L{\rm log} L
(\Omega) \to L^1(\Omega),$$ as a consequence of \eqref{G} and of a
special case of Theorem \ref{singularcor}, with $L^A(\Omega)=L{\rm
log} L (\Omega)$ and $L^B(\Omega)=L^1(\Omega)$. Thus, $\mathcal
B_\Omega f_k \to \mathcal B_\Omega f$ in $L^1(\Omega, \rn)$ and
$H_{ij}f_k \to H_{ij}f$ in $L^1(\Omega)$. This implies that $\bfu
\in W^{1,1}_0(\Omega,\R^n)$, and \eqref{partialbog} and \eqref{div}
hold.
\par\noindent
 By
Theorem \ref{singularcor}, the singular integral operator defined by
the first addend on the right-hand side of \eqref{H} is bounded from
$L^A(\Omega )$ into $L^B(\Omega )$. By inequality \eqref{G}, the
operator defined by the second addend on the right-hand side of
\eqref{H} has (at least) the same boundedness properties as a Riesz
potential operator with kernel $\frac 1{|x-y|^{n-1}}$. Such an
operator is bounded in $L^1(\Omega)$ and in $L^\infty (\Omega)$,
with norms depending only on $|\Omega|$ and on $n$. An interpolation
theorem by Calderon  \cite[Theorem 2.12, Chap. 3]{BS} then ensures
that  it is also bounded  from $L^A(\Omega)$ into $L^A(\Omega )$,
and hence, a fortiori, from $L^A(\Omega)$ into $L^B(\Omega )$, with
norm depending on $n$ and $|\Omega|$. Finally, the operator given by
the last addend on the right-hand side of \eqref{H} is pointwise
bounded (in absolute value) by $|f(x)|$. Thus, it is bounded from
$L^A(\Omega)$ into $L^A(\Omega )$, and hence  from $L^A(\Omega)$
into $L^B(\Omega )$.
Equations \eqref{bogov} and \eqref{grad} are thus established.
\par\noindent
Inequality \eqref{gradint} can be derived from \eqref{grad} via a
scaling argument analogous to that which leads to
\eqref{singularintegr} from \eqref{singularorlicz} -- see the Proof
of Theorem \ref{singularcor}.
%
%
\qed

We need a last preliminary result  in preparation for the proof of
Theorem \ref{thm:3.2}.

\begin{props}\label{supapprox} Let $\Omega$ be an open subset
in $\rn$ such that $|\Omega |< \infty$, and let $A$ be a Young
function. Assume that $u \in L^A(\Omega )$. Then:
\begin{equation}\label{supapprox1}
\sup_{v \in L^{\widetilde A}(\Omega )} \frac{\int _\Omega uv\,
dx}{\|v\|_{L^{\widetilde A}(\Omega )}} = \sup_{\varphi \in C^\infty
_0(\Omega )} \frac{\int _\Omega u\varphi\,
dx}{\|\varphi\|_{L^{\widetilde A}(\Omega )}},
\end{equation}
and
\begin{equation}\label{supapprox2}
\sup_{v \in L^{\widetilde A}_\bot(\Omega )} \frac{\int _\Omega uv\,
dx}{\|v\|_{L^{\widetilde A}(\Omega )}} = \sup_{\varphi \in C^\infty
_{0,\bot}(\Omega )} \frac{\int _\Omega u\varphi\,
dx}{\|\varphi\|_{L^{\widetilde A}(\Omega )}}.
\end{equation}
\end{props}

Note that equation \eqref{supapprox1} is well known under the
assumption that $A \in \nabla _2$ near infinity, namely $\widetilde
A \in \Delta _2$ near infinity, since $C^\infty _0 (\Omega)$ is
dense in $L^{\widetilde A}(\Omega )$ in this case. Equation
\eqref{supapprox2} also easily follows from this property when $A
\in \nabla _2$ near infinity. The novelty of Proposition
\ref{supapprox} is in the arbitrariness of $A$.

\smallskip
\par\noindent
{\bf Proof of Proposition \ref{supapprox}}. Consider first
\eqref{supapprox1}. It clearly suffices to show that
\begin{equation}\label{supapprox3}
\sup_{v \in L^{\widetilde A}(\Omega )} \frac{\int _\Omega uv\,
dx}{\|v\|_{L^{\widetilde A}(\Omega )}} = \sup_{v \in L^\infty(\Omega
)} \frac{\int _\Omega uv\, dx}{\|v\|_{L^{\widetilde A}(\Omega )}},
\end{equation}
and
\begin{equation}\label{supapprox4}
\sup_{v \in L^\infty(\Omega )} \frac{\int _\Omega uv\,
dx}{\|v\|_{L^{\widetilde A}(\Omega )}} = \sup_{\varphi \in C^\infty
_0(\Omega )} \frac{\int _\Omega u\varphi\,
dx}{\|\varphi\|_{L^{\widetilde A}(\Omega )}}.
\end{equation}
Given any $v \in L^{\widetilde A}(\Omega )$, define, for $k \in \N$,
the function $v_k : \Omega \to \R$  as
\begin{equation}\label{vk} v_k =
{\rm sign}(v)\,\min \{|v|, k\}.
\end{equation} Clearly, $v_k \in L^\infty
(\Omega)$, and $0 \leq |v_k| \nearrow |v|$ a.e. in $\Omega$ as $k
\to \infty$. Hence,
$$\int _\Omega |uv_k|\, dx \nearrow \int _\Omega |uv|\,
dx \quad \hbox{as $k \to \infty$,}$$ by the monotone convergence
theorem for integrals, and, by the Fatou property of the Luxemburg
norm,
$$\|v_k\|_{L^{\widetilde A}(\Omega )} \nearrow \|v\|_{L^{\widetilde A}(\Omega )}\quad \hbox{as $k \to
\infty$.}$$
Thus, since
$$\sup_{v \in L^{\widetilde A}(\Omega )}
\frac{\int _\Omega uv\, dx}{\|v\|_{L^{\widetilde A}(\Omega )}} =
\sup_{v \in L^{\widetilde A}(\Omega )} \frac{\int _\Omega |uv|\,
dx}{\|v\|_{L^{\widetilde A}(\Omega )}},$$
 equation \eqref{supapprox3} follows.
\par\noindent
As far as \eqref{supapprox4} is concerned, consider an increasing
sequence of compact sets $E_k$  such that ${\rm dist} (E_k, \rn
\setminus \Omega) \geq \tfrac 2k$, $E_k \subset E_{k+1} \subset
\Omega$ for $k \in \N$, and $\cup _k E_k = \Omega$.
Moreover, let $\{\varrho _k\}$ be a family of (nonnegative) smooth
mollifiers in $\rn$, such that ${\rm supp} \varrho_k \subset
B_{\frac 1k}(0)$ and $\int _{\rn} \varrho_k \, dx =1$ for $k \in \N$.
 Given $v \in
L^\infty (\Omega )$,
define $w_k : \rn \to \R$ as
$$w_k = \begin{cases} v\quad & \hbox{in $E_k$,} \\
0 & \hbox{elsewhere,}
\end{cases}$$
and $\varphi_k : \rn \to \R$ as
\begin{equation}\label{phik}
\varphi_k(x) = \int _{\rn} w_k(y)\varrho _k (x-y)\, dy \quad
\hbox{for $x \in \rn$.}
\end{equation}
 Classical properties of mollifiers
ensure that
$$\varphi_k \in C^\infty_0(\Omega), \quad \varphi_k \to v \,\,\hbox{a.e. in
$\Omega$ as $k \to \infty$}, \quad \|\varphi_k\|_{L^\infty(\Omega )}
\leq \|v\|_{L^\infty(\Omega )} \,\,\hbox{for $k \in \N$.}$$ Thus, if
$u \in L^A(\Omega)$, then
\begin{equation}\label{supapprox5}
\int _\Omega u\varphi_k\, dx \to \int _\Omega uv\, dx \quad \hbox{as
$k \to \infty$,}
 \end{equation}
 by the dominated convergence theorem
for integrals. Moreover,
\begin{equation}\label{supapprox6}
\|\varphi_k\|_{L^{\widetilde A}(\Omega )} \to \|v\|_{L^{\widetilde
A}(\Omega )}\quad \hbox{as $k \to \infty$.}
\end{equation}
Indeed, by dominated convergence and the definition of Luxemburg
norm,
\begin{equation*}
\int _\Omega \widetilde A
\bigg(\frac{|\varphi_k|}{\|v\|_{L^{\widetilde A}(\Omega )}}
\bigg)\,dx \to \int _\Omega \widetilde A
\bigg(\frac{|v|}{\|v\|_{L^{\widetilde A}(\Omega )}} \bigg)\,dx \leq
1\quad
 \hbox{as $k
\to \infty$.}
\end{equation*}
In particular, for every $\varepsilon >0$, there exists
$k_\varepsilon$ such that
$$\int _\Omega \widetilde A
\bigg(\frac{|\varphi_k|}{\|v\|_{L^{\widetilde A}(\Omega )}}
\bigg)\,dx < 1 + \varepsilon\quad \hbox{if $k >k_\varepsilon$.}$$
Hence, by the arbitrariness of $\varepsilon$ and the definition of
Luxemburg norm,
\begin{equation}\label{supapprox7}
\liminf _{k \to \infty} \|\varphi_k\|_{L^{\widetilde A}(\Omega )}
\geq \|v\|_{L^{\widetilde A}(\Omega )}.
\end{equation}
We also have that
\begin{equation}\label{supapprox8}
\limsup _{k \to \infty} \|\varphi_k\|_{L^{\widetilde A}(\Omega )}
\leq \|v\|_{L^{\widetilde A}(\Omega )}.
\end{equation}
Indeed, assume  that \eqref{supapprox8} fails. Then, there exists
$\sigma >0$ and a subsequence of $\{\varphi_k\}$, still denoted by
$\{\varphi_k\}$, such that
\begin{equation*}
1< \int _\Omega \widetilde A
\bigg(\frac{|\varphi_k|}{\|v\|_{L^{\widetilde A}(\Omega )}+ \sigma}
\bigg)\,dx \to \int _\Omega \widetilde A
\bigg(\frac{|v|}{\|v\|_{L^{\widetilde A}(\Omega )}+ \sigma}
\bigg)\,dx \leq 1\,,
\end{equation*}
a contradiction. Equation \eqref{supapprox6} follows from
\eqref{supapprox7} and \eqref{supapprox8}. Coupling
\eqref{supapprox5} with \eqref{supapprox6} yields
\eqref{supapprox4}. The proof of \eqref{supapprox1} is complete.
\par\noindent
The proof of \eqref{supapprox2} follows along the same lines, and,
in particular, via the equations
\begin{equation}\label{supapprox10}
\sup_{v \in L^{\widetilde A}_\bot(\Omega )} \frac{\int _\Omega uv\,
dx}{\|v\|_{L^{\widetilde A}(\Omega )}} = \sup_{v \in
L^\infty_\bot(\Omega )} \frac{\int _\Omega uv\,
dx}{\|v\|_{L^{\widetilde A}(\Omega )}},
\end{equation}
and
\begin{equation}\label{supapprox11}
\sup_{v \in L^\infty_\bot (\Omega )} \frac{\int _\Omega uv\,\
dx}{\|v\|_{L^{\widetilde A}(\Omega )}} = \sup_{\varphi \in C^\infty
_{0, \bot}(\Omega )} \frac{\int _\Omega u\varphi\,
dx}{\|\varphi\|_{L^{\widetilde A}(\Omega )}}.
\end{equation}
On defining, for any $v \in L^{\widetilde A}_\bot(\Omega )$,  the
sequence of functions $\{\overline v_k\} \subset L^\infty_\bot
(\Omega )$ as
$$\overline v_k = v_k - (v_k)_\Omega $$
for $k \in \N$, where $v_k$ is given by \eqref{vk}, one can prove
equation \eqref{supapprox10} via a slight variant of the argument
employed for \eqref{supapprox3}. Here, one has to use the fact that
 $(v_k)_\Omega \to
0$ as $k \to \infty$.
\par\noindent
Similarly, equation \eqref{supapprox11} can be established similarly
to \eqref{supapprox4} on replacing, for any given $v \in
L^\infty_\bot (\Omega )$, the sequence $\{\varphi _k\}$ defined  by
\eqref{phik} with
 the
sequence $\{\overline \varphi _k\} \subset C^\infty _{0,
\bot}(\Omega )$ defined as
$$\overline \varphi _k = \varphi _k - (\varphi _k)_\Omega \,\psi
\quad \hbox{for $k \in \N$,}$$
where  $\psi$ is any function in
$C^\infty _0(\Omega)$ such that $\int _\Omega \psi \, dx =1$. Note
that now, for every $\varepsilon >0$, there exists $k_\varepsilon
\in \N$ such that $\|\overline \varphi _k\|_{L^\infty (\Omega )}
\leq \|v\|_{L^\infty (\Omega )}+\varepsilon$, provided that $k
>k_\varepsilon$.
\qed

\smallskip
\par\noindent
{\bf Proof of Theorem \ref{thm:3.2}}. Let $u \in L^1(\Omega )$.
Then
\begin{align}\label{sup1}
\|u-u_\Omega\|_{L^B(\Omega)} & = \sup_{v \in L^{\widetilde A}(\Omega
)} \frac{\int _\Omega (u-u_\Omega)\,v\, dx}{\|v\|_{L^{\widetilde
B}(\Omega )}}
 =
\sup_{v \in L^{\widetilde A}(\Omega )} \frac{\int _\Omega
(u-u_\Omega)\,(v-v_\Omega)\, dx}{\|v\|_{L^{\widetilde B}(\Omega )}}
\\ \nonumber & =
\sup_{v \in L^{\widetilde A}(\Omega )} \frac{\int _\Omega
u\,(v-v_\Omega)\, dx}{\|v\|_{L^{\widetilde B}(\Omega )}}
 \leq
3\sup_{v \in L^{\widetilde A}(\Omega )} \frac{\int _\Omega
u\,(v-v_\Omega)\, dx}{\|v-v_\Omega\|_{L^{\widetilde B}(\Omega )}}
\\ \nonumber & =
3\sup_{v \in L^{\widetilde A}_\bot(\Omega )} \frac{\int _\Omega
u\,v\, dx}{\|v\|_{L^{\widetilde B}(\Omega )}}
 =
 3\sup_{\varphi \in C^\infty
_{0,\bot}(\Omega )} \frac{\int _\Omega u\varphi\,
dx}{\|\varphi\|_{L^{\widetilde A}(\Omega )}}.
\end{align}
Note that the inequality in \eqref{sup1} holds since, by the first
inequality in \eqref{AAtilde},
\begin{align*}\|v-v_\Omega\|_{L^{\widetilde B}(\Omega )}& \leq \|v\|_{L^{\widetilde B}(\Omega
)}+ \|v_\Omega\|_{L^{\widetilde B}(\Omega )} \leq
\|v\|_{L^{\widetilde B}(\Omega )}+ |v_\Omega|\|1\|_{L^{\widetilde
B}(\Omega )}
\\ \nonumber & \leq
\|v\|_{L^{\widetilde B}(\Omega )}+ \tfrac 2{|\Omega|}
\|v\|_{L^{\widetilde B}(\Omega )} \|1\|_{L^{ B}(\Omega
)}\|1\|_{L^{\widetilde B}(\Omega )} \\ \nonumber & =
\|v\|_{L^{\widetilde B}(\Omega )}+ \tfrac 2{|\Omega|}
\|v\|_{L^{\widetilde B}(\Omega )} \tfrac 1{B^{-1}(|\Omega |)} \tfrac
1{\widetilde B^{-1}(|\Omega |)} \leq 3 \|v\|_{L^{\widetilde
B}(\Omega )},
\end{align*}
and the last equality in \eqref{sup1} relies upon
\eqref{supapprox2}. By Theorem \ref{thm:3.1}, applied with $A$ and
$B$ replaced with $\widetilde B$ and $\widetilde A$, respectively,
there exists a constant $C=C(\Omega , c)$ such that
\begin{align}\label{sup2}
\sup_{\varphi \in C^\infty _{0,\bot}(\Omega )} \frac{\int _\Omega
u\varphi\, dx}{\|\varphi\|_{L^{\widetilde A}(\Omega )}} & =
 \sup_{\varphi\in L^{\widetilde{B}}_{\perp}(\Omega)} \frac{\int_\Omega u\,\Div(\mathcal B_\Omega \varphi)\,
 dx}{\|\varphi\|_{L^{\widetilde{B}}(\Omega)}}\leq\,C\, \sup_{\varphi \in L^{\widetilde{B}}_{\perp}(\Omega)}
\frac{\int_\Omega u\,\Div (\mathcal B_\Omega \varphi)\,
dx}{\|\nabla\mathcal B_\Omega
\varphi\|_{L^{\widetilde{A}}(\Omega,\R^{n\times n})}}
\\ \nonumber
&\leq\,C\, \sup_{\bfvarphi\in C^\infty_0(\Omega , \R^n)}
\frac{\int_\Omega
u\,\Div\bfvarphi\,dx}{\|\nabla\bfvarphi\|_{L^{\widetilde{A}}(\Omega,\R^{n\times
n})}}\,dx =C\,\|u\|_{W^{-1,A}(\Omega)}.
\end{align}
%
%
%
%
%
%
%
%
%
The first inequality in \eqref{negineq} follows from \eqref{sup1}
and \eqref{sup2}. The second inequality is trivial, since
\begin{align*}
\,\|u\|_{W^{-1,A}(\Omega)}& = \sup_{\bfvarphi\in C^\infty_0(\Omega ,
\R^n)} \int_\Omega
\frac{u\,\Div\bfvarphi}{\|\nabla\bfvarphi\|_{L^{\widetilde{A}}(\Omega,\R^{n\times
n})}}\,dx = \sup_{\bfvarphi\in C^\infty_0(\Omega, \R^n)}
 \frac{\int_\Omega\big(u-u_\Omega\big)\,\Div\bfvarphi \,dx}{\|\nabla\bfvarphi\|_{L^{\widetilde{A}}(\Omega,\R^{n\times n})}}\\
& \leq  C\,\sup_{\bfvarphi \in C^\infty_0(\Omega , \R^n)}
\frac{\int_\Omega \big(u-u_\Omega\big)\,\Div\bfvarphi \,
dx}{\|\Div\bfvarphi\|_{L^{\widetilde{A}}(\Omega)}}\,dx \leq
C\,\sup_{\varphi\in C^\infty_{\perp}(\Omega)}
\frac{\int_\Omega\big(u-u_\Omega\big)\,\varphi\,
dx}{\|\varphi\|_{L^{\widetilde{A}}(\Omega)}}\,dx
\\&  \leq  2 C\|u-u_\Omega\|_{L^A(\Omega)},
\end{align*}
for some constant $C=C(n)$.
\qed

\section{Nonlinear systems in fluid mechanics}\label{fluids}

In many customary mathematical models,  the stationary flow of a
homogeneous incompressible fluid in a bounded domain $\Omega\subset
\R^n$ is described by a system with the structure
\eqref{eq:1.1syst}. With a slight abuse of notation with respect to
\eqref{eq:1.1syst}, we also denote by $\bfS : \rmat \to \rmat$ the
function, acting on the symmetric gradient $\ep(\bfv)$ of
the velocity field $\bfv$, which yields the stress deviator of the
fluid. Thus, we shall consider systems of the form
\begin{align}\label{eq:4.1}
\begin{cases}
-\Div \bfS (\ep(\bfv))+\varrho\Div\big(\bfv\otimes\bfv\big)+\nabla \pi=\varrho\Div\bfF \quad & \hbox{in $\Omega$,}\\
\Div \bfv=0 & \hbox{in $\Omega$,}\\
\bfv=0 & \hbox{on $\partial \Omega$.}\end{cases}
%
%
\end{align}
In the simplest case of a Newtonian fluid, the function $\bfS$ is
linear, and   $\Div \bfS(\ep(\bfv)) = \Delta \bfv$, the Laplacian of
$\bfv$.
 However only fluids with an easy molecular structure, such
as water, oil, and several gases are governed by this low. More
complex liquids are not, and are called Non-Newtonian fluids -- see
e.g. \cite{AM, BAH}. The most common nonlinear  model among
rheologists is the   power law model, corresponding to the choice
\begin{align}\label{eq:4.2}
\bfS(\bfxi) =\nu_0\big(\kappa _0+|\bfxi|\big)^{p-2} \bfxi \quad
\hbox{for $\xi \in \R^{n\times n}$.}
\end{align}
Here,  $\nu_0 \in (0, \infty)$ and $\kappa _0 \in [0, \infty)$ are
constants, and $p\in(1,\infty)$ is an exponent which need to be
specified via physical experiments. An extensive list of specific
$p$-values for different fluids can be found in \cite{BAH}.
\par\noindent
%
A more general constitutive law for Non-Newtonian fluids, which
allows for non-polynomial type nonlinearities,  takes the form
\begin{align}\label{eq:4.4''}
\bfS(\bfxi) =\frac{\Phi '(|\bfxi|)}{|\bfxi|} \bfxi \quad \hbox{for $\bfxi \in
\R^{n\times n}$}\,,
\end{align}
where $\Phi$ is a Young function.
\par\noindent
In various instances of interest in applications, the term
\begin{equation}\label{vv}
\Div\big(\bfv\otimes\bfv\big)=(\nabla\bfv)\bfv
\end{equation}
is negligible in \eqref{eq:4.1}, compared with the other terms
appearing in the first equation. This is the case, for example, if
the modulus of the velocity $\bfv$ is small. Another situation where
the role of the term \eqref{vv} is immaterial is that of plastic or
pseudo-plastic fluids. Indeed, \eqref{vv} accounts for the inner
rotation in the fluid flow, and for such fluids the impact of this
term is very limited. Dropping the term \eqref{vv} reduces
\eqref{eq:4.1} to the simplified system
\begin{align}\label{eq:4.1'}
\begin{cases}
-\Div \bfS (\ep(\bfv))+\nabla \pi=\varrho\Div\bfF &  \quad \hbox{in $\Omega$,}\\
\Div \bfv=0 &\hbox{in $\Omega$,}\\
\bfv=0 & \hbox{on $\partial \Omega$.}\end{cases}
%
%
%
%
\end{align}
A standard approach to (\ref{eq:4.1}) or (\ref{eq:4.1'}) consists in
two steps. Firstly,  a velocity field $\bfv$ is exhibited  such that
\begin{align}\label{Hphi}
\int_\Omega \bfH:\nabla\bfphi\,dx&=0
\end{align}
for every $\bfvarphi\in C^{\infty}_{0,\Div}(\Omega, \rn)$, where
either
\begin{equation}\label{Hv}
\bfH=\bfS (\ep(\bfv)) +\rho\,\bfF\,,
\end{equation}
or
\begin{align}\label{Hred}
\bfH=\bfS (\ep(\bfv)) +\rho\,\bfF\,-\rho\,\bfv\otimes\bfv\,,
\end{align}
according to weather the convective term $\bfv\otimes\bfv$ is
included in the model or not. Here, \lq\lq \,: " stands for scalar
product between matrices, and $C^{\infty}_{0,\Div}(\Omega, \rn)$
denotes the space of compactly supported, infinitely differentiable
$\rn$-valued functions whose divergence vanishes in $\Omega$. The
function $\bfv$ belongs to a proper Sobolev type space depending on
the constitutive law underlying the definition of the function
$\bfS$.
 Secondly,  the pressure $\pi$ is
reconstructed.
 \par
 A discussion of  the first issue falls beyond the scopes of the present paper, and will not be addressed
 here. Let us just mention that the standard power type model \eqref{eq:4.2}
 has been investigated in the classical contributions  \cite{La, La2, La3, Li}, and in the recent papers \cite{FMS1, FMS2, DMS, DRW, BrDS}.
 Stationary
  flows of fluids whose constitutive law satisfies (\ref{eq:4.4''}) with a Young-function $\Phi\in\Delta_2\cap\nabla_2$ are studied
  in \cite{BrF, Br, DK}. An unconventional
constitutive law, where $\bfS$ has the form \eqref{eq:4.4''} with
$\Phi ' (t) \approx \log (1+t)$ near infinity, and hence $\Phi
\notin \nabla _2$ near infinity, was introduced by Eyring in
\cite{E}, where it is assumed that
\begin{align}\label{eq:4.3}
\bfS(\bfxi)= \nu_0 \frac{\mathrm{arsinh} (\lambda_0 |\bfxi|)}{\lambda_0
|\bfxi|}\bfxi \quad \hbox{for $\bfxi \in \R ^{n\times n}$,}
\end{align}
for some physical constants $\nu_0,\lambda>0$. Similar results are
due to Prandtl (see e.g. \cite{BB} for an overview on this kind of
models). An analysis of the simplified system \eqref{eq:4.1'} for
the Eyring-Prandtl model is the object of \cite{FuS}, whereas the
complete system (\ref{eq:4.1}), in the  case $n=2$, is considered in
\cite{BrDF}.
\par
In the remaining part of this paper, we focus, instead, on the
second question, namely the reconstruction of the pressure $\pi$ in
a correct Orlicz space. In case of fluids governed by a general
constitutive low of the form \eqref{eq:4.4''}, the function $\bfH$
belongs to some Orlicz space
 $L^A(\Omega,\R^{n\times n})$. If   $A\in \Delta_2 \cap \nabla_2$, then $\pi \in L^A(\Omega)$ as well. However,
in general, one can only expect that $\pi$ belongs to some larger
Orlicz space $L^B(\Omega)$. The balance  between the Young functions
$A$ and $B$ is determined by conditions \eqref{1.1} and \eqref{1.2},
as stated in the following result.

\begin{theorem}
\label{thm:4.1} 
 Let $A$ and $B$ be Young functions fulfilling
\eqref{1.1} and \eqref{1.2}. Let $\Omega$ be a bounded domain with
the cone property in $\R ^n$, $n \geq 2$.
 Assume that $\bfH\in L^A(\Omega,\R^{n\times n})$
and satisfies
\begin{align*}
\int_\Omega \bfH:\nabla\bfphi\,dx&=0
\end{align*}
for every $\bfvarphi\in C^{\infty}_{0,\Div}(\Omega, \rn)$. Then
there exists a unique function $\pi\in L_\bot^B(\Omega)$ such that
\begin{align}\label{eq:thm4.1}
\int_\Omega \bfH:\nabla\bfphi\,dx&=\int_\Omega
\pi\,\Div\bfvarphi\,dx
\end{align}
for every $\bfvarphi\in C^{\infty}_{0}(\Omega,\R^n)$. Moreover,
there exists a constant $C=C(\Omega, c)$ such that
\begin{align}\label{7}
\|\pi\|_{L^B(\Omega)}\leq  C\|\bfH-\bfH
_\Omega\|_{L^A(\Omega,\R^{n\times n})},
\end{align}
and
\begin{align}\label{8}
\int _\Omega B(|\pi|)\, dx \leq \int _\Omega A(C|\bfH- \bfH
_\Omega|)\,dx.
\end{align}
Here, $c$ denotes the constant appearing in \eqref{1.1} and
\eqref{1.2}.
\end{theorem}

In particular,  Theorem \ref{thm:4.1} reproduces, within a  unified
framework, various results appearing in the literature.  For
instance, when the constitutive relation \eqref{eq:4.2} is in force,
 the function $A(t)$ is just a power $t^q$, where the exponent $q >1$,
 and depends on $p$, on $\bfF$, and on whether the system \eqref{eq:4.1} or
 \eqref{eq:4.1'} is taken into account. In any case,   $L^A(\Omega, \R
^{n\times n})$ agrees with the Lebesgue space $L^q (\Omega, \R
^{n\times n})$, and Theorem \ref{thm:4.1} recovers the fact that
$\pi$ belongs to the same Lebesgue space $L^q(\Omega)$.
\par\noindent As far as the simplified system \eqref{eq:4.1'} for the
Eyring-Prandtl model \eqref{eq:4.3} is concerned, under appropriate
assumptions on $\bfF$ one has that $\bfH \in {\rm exp}L(\Omega,
\R^{n\times n})$. Hence, via Theorem \ref{thm:4.1}, we infer the
existence of a pressure   $\pi\in {\rm exp}L^{\frac 12}(\Omega)$.
More generally, if  $\bfH \in {\rm exp}L^\beta(\Omega, \R^{n\times
n})$ for some $\beta >0$, one has that $\pi\in {\rm
exp}L^{\beta/(\beta+1)}(\Omega)$. The complete system \eqref{eq:4.1}
for the Eyring-Prandtl model, in the $2$-dimensional case, admits a
weak solution  $\bfv$ such that $\bfv\otimes\bfv\in L {\rm log}
L^2(\Omega,\R^{n\times n})$ and hence  $\bfH\in  L {\rm log}
L^2(\Omega, \R^{n\times n})$ \cite{BrDF}. Again, one cannot expect
that the pressure $\pi$ belongs to  the same space. In fact, Theorem
\ref{thm:4.1} yields the existence of a pressure  $\pi\in
 L {\rm
log} L(\Omega)$, thus reproducing a result from \cite{BrDF}. In
general, if $\bfH\in L {\rm log} L^\alpha (\Omega, \R^{n\times n})$
for some $\alpha \geq 1$, then we obtain that $\pi\in
 L {\rm log} L^{\alpha -1}(\Omega)$.

\smallskip
\par\noindent
{\bf Proof of Theorem \ref{thm:4.1}}. By De Rahms Theorem, in the
version of \cite{Si},
there exists a distribution $\Xi$ such that
\begin{align}\label{dera}
\int_\Omega \bfH:\nabla\bfphi\,dx&=\Xi(\Div\bfvarphi)
\end{align}
for every $\bfvarphi\in C^{\infty}_{0}(\Omega, \rn)$. Replacing
$\bfvarphi$ with $\mathcal B_\Omega\big(\varphi- \varphi
_\Omega\big)$ in \eqref{dera}, where $\varphi \in C^\infty _0(\Omega
)$, yields
\begin{align*}
\int_\Omega \bfH:\nabla\mathcal B_\Omega\big(\varphi- \varphi
_\Omega\big)\,dx&=\Xi\big(\varphi- \varphi _\Omega\big)
\end{align*}
for every $\varphi\in C^{\infty}_{0}(\Omega)$. We claim that the
linear functional $C^{\infty}_{0}(\Omega) \ni
\varphi\mapsto\Xi\big(\varphi- \varphi_\Omega\big)$ is bounded on
$C^{\infty}_{0}(\Omega)$ equipped with the $L^\infty (\Omega)$ norm.
Indeed,  by (\ref{1.1}), one has that $L^A(\Omega, \rmat)\rightarrow
L\mathrm{log}L(\Omega, \rmat)$.
Moreover, by a special case of Theorem \ref{thm:3.1}, $\nabla
\mathcal B_\Omega : L_\bot ^\infty (\Omega) \to
\mathrm{exp}L(\Omega, \rmat)$. Thus, since $L\mathrm{log}L(\Omega,
\rmat)$ and $\mathrm{exp}L(\Omega, \rmat)$ are Orlicz spaces built
upon Young functions which are conjugate of each other,
\begin{align}\label{march1}
\bigg|\int_\Omega \bfH:\nabla\mathcal B_\Omega\big(\varphi - \varphi
_\Omega\big)\,dx\bigg| & \leq C \|\bfH \|_{L{\rm log}L(\Omega,
\rmat)} \|\nabla\mathcal B_\Omega\big(\varphi - \varphi
)_\Omega\big)\|_{\mathrm{exp}L(\Omega, \rmat)}
\\ \nonumber & \leq C'  \|\bfH \|_{L^A(\Omega, \rmat)} \|\varphi - \varphi
_\Omega\|_{L^\infty (\Omega)} \\ \nonumber & \leq C''  \|\bfH
\|_{L^A(\Omega, \rmat)} \|\varphi\|_{L^\infty (\Omega)}\,,
\end{align}
for every $\varphi\in C^{\infty}_{0}(\Omega)$, where $C=C(|\Omega|,
n)$ and $C'=C'(\Omega, c)$. Hence, the relevant functional can be
continued to a bounded linear functional on $\varphi\in
C^0_{0}(\Omega)$, with the same norm.
\par\noindent
Now, as a consequence of Riesz's representation Theorem, there
exists  a Radon measure $\Xi$ such that
\begin{align*}
\Xi\big(\varphi- \varphi _\Omega\big)=\int_\Omega \varphi\,d\mu
\end{align*}
for every $\varphi\in C^0_{0}(\Omega)$. Fix any open set $E\subset
\Omega$. By Theorem \ref{thm:3.1} again, there exists a constant $C$
such that
\begin{align}\label{ac}
\mu(E)&=\sup_{\varphi\in
C^0_{0}(E),\,\,\|\varphi\|_\infty=1}\Xi\big(\varphi-
\varphi_\Omega\big) =\sup_{\varphi\in
C^0_{0}(E),\,\,\|\varphi\|_\infty=1}\int_\Omega
\bfH:\nabla\mathcal B_\Omega\big(\varphi- \varphi_\Omega\big)\,dx\\
\nonumber &\leq \sup_{\varphi\in
C^0_{0}(E),\,\,\|\varphi\|_\infty=1}\|\bfH\|_{L{\rm log}L(E,\R^{n\times n})}
\|\nabla\mathcal
B_\Omega\big(\varphi-\varphi_\Omega\big)\|_{\mathrm{exp}L(\Omega,\R^{n\times n})}\\
\nonumber &\leq C\,\sup_{\varphi\in
C^0_{0}(E),\,\,\|\varphi\|_\infty=1}\|\bfH\|_{L{\rm
log}L(E,\R^{n\times n})}\|\varphi-(\varphi)_\Omega\|_{L^{\infty}(\Omega)} \leq
C\,\|\bfH\|_{L{\rm log}L(E,\R^{n\times n})}.
\end{align}
One can verify that the norm $\|\cdot \|_{L{\rm log}L(E)}$ is
absolutely continuous, in the sense that for every $\varepsilon >0$
there exists $\delta
>0$ such that $\|\bfH\|_{L{\rm log}L(E,\R^{n\times n})} <
\varepsilon$ if $|E| < \delta$, and since any Lebesgue measurable
set can be approximated from outside by open sets, inequality
\eqref{ac} implies that the measure $\mu$ is absolutely continuous
with respect to the Lebesgue measure. Hence,
 $\mu$ has a density with respect to the Lebesgue measure. So
$\Xi$ can be represented by a function  $\pi \in L^1(\Omega)$
fulfilling (\ref{eq:thm4.1}) holds. The function $\pi$ is uniquely
determined if we assume that $\pi _\Omega=0$.
 By this assumption, Theorem
\ref{thm:3.2}, and equation (\ref{eq:thm4.1}) we have that
\begin{align*}
\|\pi\|_{L^B(\Omega)}&\leq C\|\nabla\pi\|_{W^{-1,A}(\Omega, \rn)} =
C\sup_{\bfvarphi\in C^\infty_0(\Omega, \rn)}
\frac{\int_\Omega\pi\,{\rm div}\, \bfvarphi\,dx}{\|\nabla
\bfvarphi\|_{L^{\widetilde{A}}(\Omega, \rmat)}}
\\
&= C\sup_{\bfvarphi\in C^\infty_0(\Omega, \rn)} \frac{\int_\Omega
\bfH \, : \, \nabla \bfvarphi\, dx}{\|\nabla
\bfvarphi\|_{L^{\widetilde{A}}(\Omega, \rmat)}} =C\sup_{\bfvarphi\in
C^\infty_0(\Omega, \rn)} \frac{\int_\Omega(\bfH- \bfH _\Omega)\, :
\, \nabla \bfvarphi\, dx}{\|\nabla
\bfvarphi\|_{L^{\widetilde{A}}(\Omega, \rmat)}}
\\ & \leq 2 C
\|\bfH- \bfH _\Omega\|_{L^A(\Omega, \rmat)},
%
\end{align*}
where $C=C(\Omega , c)$. This proves inequality \eqref{7}.
Inequality \eqref{8} follows from \eqref{7}, on replacing $A$ and
$B$ with $k A$ and $k B$, respectively, with $k = \frac 1{\int
_\Omega A(|\bfH- \bfH_\Omega|)\,dx}$, via an argument analogous to
that of the proof of \eqref{singularintegr}. \qed

\smallskip
\par
Let us turn to a further consequence of Theorem \ref{thm:4.1}, which
is related to a numerical analysis of problem \eqref{eq:4.1}, or of
its simplified version \eqref{eq:4.1'}.  We shall adopt the scheme
 of the
 finite element method   for the $p$-Stokes system exploited in
\cite{BeBeDiRu} in the special case when $\bfS$ is  given by
\eqref{eq:4.2}. In what follows, we assume that $\Omega$ is a
polyhedron.
 The
goal is to compute an approximate solution to problem \eqref{eq:4.1'}
via discretization.
%
%
To this purpose,  one needs a triangulation $\mathscr T_h$ of
$\Omega$ into  simplices   of diameter bounded by $h>0$. Recall that
a simplex in $\rn$ is the convex hull of $n+1$ points which do not
lie on the same hyperplane. We also need that such a triangulation
is regular enough for the Lipschitz constant of the functions, which
locally represent the boundaries of the relevant simplices, to be
uniformly bounded in $h$. We denote by $\mathscr P^{0, h}(\Omega)$
the space of those functions in $\Omega$ whose restriction to each
simplex of $\mathscr T_h$ is constant, and by $\mathscr P^{0,
h}_\bot(\Omega)$ its subspace of
 those functions from $\mathscr P^{0,
h}(\Omega)$ whose mean value over $\Omega$ is zero. Also,  for $\ell
\geq1$, we denote by $\mathscr P^{\ell, h}(\Omega)$
 the space of first-order weakly differentiable functions in $\Omega$ whose restriction  to each simplex of $\mathscr T_h$
 is a polynomial of degree $\ell$, and  by $\mathscr P^{\ell, h}_0(\Omega)$ its subspace of those functions from  $\mathscr P^{\ell, h}(\Omega)$
which vanish on $\partial \Omega$. Finally $\mathscr P^{\ell,
h}_\bot(\Omega)$ stands for the space of all functions from
$\mathscr P^{\ell, h}(\Omega)$ whose mean value over $\Omega$ is
zero.
 Clearly, given any  Young function $A$, one has that
$$\mathscr P^{\ell, h}(\Omega) \subset L^A(\Omega), \quad \mathscr P^{\ell,
h}_\bot(\Omega) \subset L^A_\bot (\Omega) \quad \hbox{if $\ell \geq
0$,}$$ and
$$ \mathscr P^{\ell, h}(\Omega) \subset W^{1,A}(\Omega), \quad
 \mathscr P^{\ell, h}_0(\Omega) \subset W^{1,A}_0(\Omega) \quad \hbox{if $\ell \geq
1$.}$$
  The spaces $\mathscr P^{\ell, h}(\Omega,
 \rn)$, $\mathscr P^{\ell, h}_0(\Omega, \rn)$ and $\mathscr P^{\ell,
h}_\bot(\Omega, \rn)$  of $\rn$-valued functions are defined
accordingly.
\par An important tool for the numerical analysis of problem \eqref{eq:4.1} is
a projection operator $\Pi^h$ which, for given $m \in \mathbb N \cup
\{0\}$ and $k \in \mathbb N$, is such that
%
$\Pi^h:
W^{1,1}(\Omega,
\rn)\rightarrow \mathscr P^{k, h}(\Omega, \rn)$, and enjoys the following properties:\\
 (i) $\Pi^h$ preserves zero boundary values, i.e.
\begin{equation}\label{itm:A1-2}
\Pi^h:
W^{1,1}_0(\Omega, \rn) \to  \mathscr P^{k, h}_0(\Omega, \rn);
\end{equation}
(ii) $\Pi^h$ is divergence preserving, namely,
\begin{align}\label{itm:A1}
\int_\Omega p_h\,\Div\bfu\,dx=\int_\Omega p_h\,\Div
\Pi^h\bfu\,dx\quad \hbox{for every}\, \bfu\in W^{1,1}(\Omega, \rn)
\, \hbox{and for every}\, p_h \in \mathscr P^{m, h}(\Omega);
\end{align}
(iii)  $\Pi^h$ is continuous in the $W^{1,1}$-sense, i.e. there
exists a constant $C=C(\Omega)$ such that
\begin{align}\label{itm:A2}
\dashint_{\mathcal S}|\Pi^h\bfu|\,dx+\dashint_{\mathcal S}h_\mathcal
S|\nabla \Pi^h\bfu|\,dx\leq C\,\dashint_{M_{\mathcal
S}}|\bfu|\,dx+C\,\dashint_{M_\mathcal S}h_\mathcal S|\nabla
\bfu|\,dx
\end{align}
for every $\mathcal S\subset\mathscr T_h$, and every $\bfu \in
W^{1,1}(\Omega , \rn)$. Here $h_\mathcal S$ denotes the diameter of
$\mathcal S$, and $M_{\mathcal S}$ the union of $\mathcal S$ and all
its direct neighbors in the triangulation $\mathscr T_h$.
The existence of such an operator for suitable couples
of $m$ and $k$ is standard -- see e.g. \cite[appendix]{BeBeDiRu}.
For instance, the choice  $k=2$ and $m=0$ is admissible,  whereas
$k=1$ and $m=0$ is not.

%
%
%
Now, the discrete version of problem \eqref{eq:4.1'} amounts to
finding a couple $(\bfv_h,\pi_h)\in \mathscr P^{k, h}_0(\Omega,
\rn)\times \mathscr P^{m, h}_\bot(\Omega)$ such that
\begin{align*}
\int_\Omega \bfS(\ep(\bfv_h)):\nabla \bfphi_h\,dx&=\int_\Omega
\pi_h\,\Div\bfphi_h\,dx-\rho\int_\Omega\bfF :
\nabla\bfphi_h\,dx,\end{align*} and
\begin{align*}
 \int_\Omega p_h\,\Div \bfv_h\,dx&=0,
\end{align*}
for every $(\bfvarphi_h,p_h)\in \mathscr P^{k, h}_0(\Omega,
\rn)\times \mathscr P^{m, h}_\bot(\Omega)$.
\par
This can again be accomplished in two steps. First, on setting
$$\mathscr P^{k, h}_{0, {\Div}}(\Omega, \rn)= \left\{\bfu_h\in \mathscr P^{k, h}_0(\Omega, \rn):\,\,\int_\Omega
p_h\,\Div\bfu_h\,dx=0\,\,\hbox{for every}\,\, p_h\in \mathscr P^{m,
h}(\Omega) \right\},$$ and
$$\bfH_h=\bfS(\ep(\bfv_h))+\rho\bfF\,,$$
 one has to find
$\bfv_h\in \mathscr P^{k, h}_{0, {\Div}}(\Omega, \rn)$ such that
\begin{align}\label{discrpi}
\int_\Omega \bfH_h:\nabla \bfphi_h\,dx =0
%
\quad \hbox{for every $\bfvarphi_h\in \mathscr P^{k, h}_{0,
{\Div}}(\Omega, \rn)$}.
\end{align}
The existence of a unique function $\bfv_h$ satisfying
\eqref{discrpi}  follows easily from the theory of monotone
operators. If $\bfS$ has variational structure one may equivalently
solve the corresponding strictly convex minimizing problem on the
finite dimensional function space $\mathscr P^{k, h}_{0,
\Div}(\Omega, \rn)$. Next, the pressure has be to reconstructed.
Precisely, one has to find $\pi_h\in \mathscr P^{m, h}_\bot(\Omega)$
such that
 \begin{align*}
\int_\Omega \bfH_h:\nabla\bfphi_h\,dx&=\int_\Omega \pi_h\,\Div\bfphi_h\,dx\quad \hbox{for every $\bfvarphi_h\in \mathscr P^{k, h}_0(\Omega, \rn)$.}
\end{align*}
This is the
discrete analogone of the problem from Theorem \ref{thm:4.1}.
 As far as existence is concerned, we only need to solve an algebraic linear system. The following questions then arise:
\par\noindent (i) Is the pressure $\pi_h$ unique?
\par\noindent (ii)  Does the pressure $\pi_h$ depend continuously on the data of the
problem, namely  on $\bfF$?
\par\noindent (iii) Does the family of discretized pressure functions $\{\pi_h\}$ converge when $h\rightarrow0$?
\par
The answer to all this questions follows from the so-called inf-sup
condition. Such a condition, in the standard case when $\bfS$ has a
power type growth as in \eqref{eq:4.2}, reads as
\begin{align}\label{infsupp}
\inf_{p_h\in P^{m,
h}_\bot(\Omega)}\sup_{\bfvarphi_h\in \mathscr P^{k,
h}_0(\Omega, \rn)}\frac{\int_\Omega
p_h\,\Div\bfphi_h\,dx}{\|p_h\|_{L^{p'}(\Omega)}\|\nabla
\bfvarphi_h\|_{L^p(\Omega, \rmat)}}\geq C,
\end{align}
for some positive constant $C$ independent of $h$. In fact, the
uniqueness of $\pi_h$ does not even require $C$ to be independent of
$h$.
\par
 Our next result
provides us with an Orlicz space version of \eqref{infsupp}.

\begin{theorem}
\label{thm:4.2}
 Let $A$ and $B$ be Young functions fulfilling
\eqref{1.1} and \eqref{1.2}. Assume that $\Omega$ is a polyhedron in
$\R ^n$, $n \geq 2$, and that $m \in \mathbb N \cup
\{0\}$ and $k \in \mathbb N$ are such that there exists an operator
$\Pi^h: W^{1,1}(\Omega, \rn)\rightarrow \mathscr P^{k, h}(\Omega,
\rn)$  satisfying (\ref{itm:A1-2})-(\ref{itm:A2}). Then,
\begin{align}\label{infsuporlicz}
\inf_{p_h\in \mathscr P^{m, h}_\bot(\Omega)}\sup_{\bfvarphi_h\in
\mathscr P^{k, h}_0(\Omega, \rn)}\frac{\int_\Omega
p_h\,\Div\bfphi_h\,dx}{\|p_h\|_{L^B(\Omega)}\|\nabla
\bfvarphi_h\|_{L^{\tilde{A}}(\Omega, \rmat)}} \geq C,
\end{align}
for some positive constant $C=C(\Omega, c)$, where $c$ is the
constant appearing in \eqref{1.1} and \eqref{1.2}
\end{theorem}

\smallskip
\par\noindent
{\bf Proof}. As a first step, we show that the
operator $\Pi^h$ is continuous in every Orlicz
 space, in the sense that there exists a constant $C=C(\Omega)$ such that, for
 every Young function $A$,
\begin{align}\label{itm:A2new}
\|\nabla \Pi^h\bfu\|_{L^A(\Omega, \rmat)}\leq C
\|\nabla\bfu\|_{L^A(\Omega, \rmat)}\quad \hbox{for every}\,  \bfu\in
W^{1,A}(\Omega, \rn).
\end{align}
Inequality \eqref{itm:A2new} has been established
 in \cite[Thm.
4.5]{DiRu} and in \cite[Thm. 3.2]{BeBeDiRu}
%
under the additional assumption that $A$ fulfils a global
$\Delta_2$-condition. A  variant of those proofs
shows that, in  fact, this assumption can be dropped. We outline the
argument hereafter. Since $|\nabla \Pi^h\bfu|$
belongs to a finite dimensional function space, it follows from
(\ref{itm:A2}) that
\begin{align*}
\dashint_{\mathcal S}A(h_\mathcal S|\nabla \Pi^h\bfu|)\,dx& \textcolor{blue}{\leq} \dashint_\mathcal S A\bigg(\textcolor{blue}{C}\dashint_{\mathcal S}h_\mathcal S|\nabla \Pi^h\bfu|\,dy\bigg)\,dx\\
&\leq \dashint_\mathcal S A\bigg(C'\,\dashint_{M_\mathcal
S}|\bfu|\,dy+C'\,\dashint_{M_{\mathcal S}}h_\mathcal S|\nabla
\bfu|\,dy\bigg)\,dx
\end{align*}
for some constants $C=C(\Omega)$ and
$C'=C'(\Omega)$ , and for every $\bfu \in W^{1,1}(\Omega , \rn)$.
Hence, owing to Jensen's inequality and the convexity of $A$,
\begin{align}\label{eq:stab}
\dashint_{\mathcal S}A(h_\mathcal S|\nabla \Pi^h\bfu|)\,dx &\leq
\dashint_{M_{\mathcal S}}A(C'|\bfu|)\,dy+\,\dashint_{M_\mathcal
S}A(C'h_\mathcal S|\nabla \bfu|)\,dy
\end{align}
for every  $\bfu \in W^{1,A}(\Omega , \rn)$.
%
Now, given any $\bfq\in\R^n$, we deduce from the convexity of $A$,
\eqref{itm:A1-2} and \eqref{eq:stab}
\begin{align}
\dashint_{\mathcal S}A(h_\mathcal S|\nabla\bfu-\nabla
\Pi^h\bfu|)\,dx& \leq \frac12\dashint_{\mathcal S}A(2h_\mathcal
S|\nabla(\bfu-\bfq) |)\,dx+
\frac 12\dashint_{\mathcal S}A(2h_\mathcal S|\nabla \Pi^h(\bfu-\bfq)|)\,dx\nonumber\\
&\leq \dashint_{M_{\mathcal S}}A(Ch_\mathcal S|\nabla(\bfu-\bfq)
|)\,dx+\dashint_{M_{\mathcal S}}A(C|\bfu-\bfq|)\,dx,\label{*}
\end{align}
for some constant $C=C(\Omega)$ and for every $\bfu \in
W^{1,A}(\Omega , \rn)$. Finally we choose $\bfq=\bfu_{M_{\mathcal
S}}$. By inequality \eqref{poincare} and a scaling argument as in
the proof of Theorem \ref{thm:3.1} one can show that there exists a
constant $C=C(\Omega )$ such that
\begin{equation}\label{poincareint}
\dashint_{M_{\mathcal S}}A(|\bfu-\bfq|)\,dx \leq
\dashint_{M_{\mathcal S}}A(Ch_\mathcal S|\nabla\bfu |)\,dx
\end{equation}
 for every $\bfu \in W^{1,A}(\Omega , \rn)$. Inequalities \eqref{*} and \eqref{poincareint} imply that
\begin{align}\label{eq:appr}
\dashint_{\mathcal S}A(h_\mathcal S|\nabla\bfu-\nabla
\Pi^h\bfu|)\,dx &\leq \dashint_{M_{\mathcal S}}A(C h_\mathcal
S|\nabla\bfu |)\,dx
\end{align}
for every $\bfu \in W^{1,A}(\Omega , \rn)$.
The convexity of $A$ and inequality (\ref{eq:appr}) yield
\begin{align*}
\dashint_{\mathcal S}A(h_\mathcal S|\nabla \Pi^h\bfu|)\,dx &\leq
\frac 12 \dashint_{\mathcal S}A(2h_\mathcal S|\nabla \bfu|)\,dx +
\frac 12 \dashint_{\mathcal S}A(2h_\mathcal S|\nabla \bfu - \nabla
\Pi^h\bfu|)\,dx
\\ &  \leq
\dashint_{M_{\mathcal S}}A(Ch_\mathcal S|\nabla\bfu |)\,dx
\end{align*}
for some constant $C=C(\Omega)$, and for every $\bfu \in
W^{1,A}(\Omega , \rn)$. This implies \eqref{itm:A2}, on  replacing
$\bfu$ with $h_\mathcal S^{-1}\bfu$ and summing up over all
simplexes $\mathcal S$. Note, in this connection, that the number of
neighbours of each $\mathcal S$ only depends on $n$.
\par\noindent
We are now in a position to prove inequality \eqref{infsuporlicz}.
By Theorem \ref{thm:3.2}, there exists a constant $C=C(\Omega,c)$
such that
\begin{align}\label{**}
\|p_h\|_{L^{B}(\Omega)}&\leq C\|\nabla p_h\|_{W^{-1,A}(\Omega,\rn)}
=C\sup_{\bfvarphi\in C^\infty
_0(\Omega, \rn)}\frac{\int_\Omega
p_h\,\Div\bfphi\,dx}{\|\nabla\bfvarphi\|_{L^{\widetilde{A}}(\Omega)}}
\end{align}
for every $p_h\in \mathscr P^{0, h}_\bot(\Omega)$. Owing to
\eqref{**} and to  properties  \eqref{itm:A1} and \eqref{itm:A2} of the
 operator $\Pi^h$,
\begin{align*}
\|p_h\|_{L^{B}(\Omega)}& \leq  C\sup_{\bfvarphi\in
C^\infty_0(\Omega, \rn)} \frac{\int_\Omega
p_h\,\Div\Pi^h\bfphi\,dx}{\|\nabla\bfvarphi\|_{L^{\widetilde{A}}(\Omega,\mathbb R^{n\times n})}}
\leq C'\sup_{\bfvarphi\in C^\infty_0(\Omega, \rn)}\frac{\int_\Omega p_h\,\Div \Pi^h\bfphi\,dx}{\|\nabla \Pi^h\bfvarphi\|_{L^{\widetilde{A}}(\Omega,\mathbb R^{n\times n})}}\\
&\leq C'\sup_{\bfvarphi_h\in  \mathscr P^{k, h}_0(\Omega,
\rn)}\frac{\int_\Omega
p_h\,\Div\bfphi_h\,dx}{\|\nabla\bfvarphi_h\|_{L^{\widetilde{A}}(\Omega,\mathbb R^{n\times n})}},
\end{align*}
for some constants $C=C(\Omega , c)$ and $C'=C'(\Omega , c)$, and
for every $p_h \in \mathscr P^{m, h}_\bot(\Omega)$. Hence,
\eqref{infsuporlicz} follows. \qed

\begin{remark} {\rm Since the functions $p_h$ belong to finite dimensional
spaces it would be possible to replace  $\|p_h\|_{L^{B}(\Omega)}$ in
Theorem \ref{thm:4.2} by $\|p_h\|_{L^{A}(\Omega)}$. However, the
constant then depends on the dimension of $\mathscr
P^{m,h}_\bot(\Omega)$ and hence on $h$. In that form, the result
would be of no use in the study of the convergence of the finite
element method approximation.}
\end{remark}

\begin{corollary}
\label{cor:4.1}
Let $A$, $B$ and $\Omega$ be as in Theorem \ref{thm:4.2}. Assume
that the function $\bfH_h\in L^A(\Omega, \rmat)$ fulfils
\begin{align}
\label{cor:4.1.0}\int_\Omega \bfH_h:\nabla\bfphi_h\,dx&=0
\end{align}
for every $\bfvarphi_h\in \mathscr P^{k, h}_{0,\Div}(\Omega, \rn)$.
Then:
\begin{itemize}
\item[i)] There exists a unique function $\pi_h\in \mathscr P^{m, h}_\bot(\Omega)$ such that
 \begin{align}\label{march1bis}
\int_\Omega \bfH_h:\nabla\bfphi_h\,dx&=\int_\Omega
\pi_h\,\Div\bfphi_h\,dx
\end{align}
for every $\bfvarphi_h\in \mathscr P^{k, h}_{0}(\Omega, \rn)$.
\item[ii)] There exists a constant $C$, independent of $h$, such that
\begin{align}\label{oct1}
\|\pi_h\|_{L^B(\Omega)}\leq C\|\bfH_h\|_{L^A(\Omega, \rmat)}.
\end{align}
 \item[iii)] There exists a  constant $C$, independent of $h$, such that, if $\bfH\in L^A(\Omega, \rmat)$ and $\pi\in L^A_\bot(\Omega)$ satisfy
 \begin{align}\label{march2}
\int_\Omega \bfH:\nabla\bfphi\,dx&=\int_\Omega \pi\,\Div\bfphi\,dx
\end{align}
for every $\bfvarphi\in C^\infty_0(\Omega, \rn)$,
then
\begin{align}\label{oct2}
\|\pi_h-\pi\|_{L^B(\Omega)}\leq C\Big(\|\bfH_h-\bfH\|_{L^A(\Omega,
\rmat)}+\inf_{\mu_h\in\mathscr P^{0,
h}_\bot(\Omega)}\|\mu_h-\pi\|_{L^A(\Omega)}\Big).
\end{align}
\end{itemize}
\end{corollary}

\smallskip
\par\noindent
 {\bf Proof}. (i) Consider a basis $\{p_h^j\}$, $i=1, \dots ,
N_h$ of $\mathscr P^{m, h}_\bot(\Omega)$ and a basis
$\{\bfvarphi_h^j\}$, $j=1, \dots , M_h$ of $\mathscr P^{k,
h}_0(\Omega, \rn)$. Note that $N_h \leq M_h$. Then the problem
 \begin{align*}
\int_\Omega \bfH_h:\nabla\bfphi_h\,dx&=\int_\Omega \pi_h\,\Div\bfphi_h\,dx
\end{align*}
for every $\bfvarphi_h\in \mathscr P^{k, h}_{0}(\Omega, \rn)$ is
equivalent to the algebraic linear system $\mathcal A \bfz=\bfb$,
where
\begin{align*}
\mathcal A_{ij}&=\int_\Omega p^j_h\,\Div\bfphi^i_h\,dx,\quad
b_i=\int_\Omega \bfH_h:\nabla\bfphi^i_h\,dx.
\end{align*}
Such system has a solution,  provided that $\bfb$ belongs to the
image of the matrix $\mathcal A$, or, equivalently, if it is
orthogonal to  $\Ker(\mathcal A^T)$. Now, observe that a vector
$\bfy$ belongs to $\Ker(\mathcal A^T)$ if and only if
\begin{align*}
0=\sum_i\mathcal A_{ij}y_i=\int_\Omega p^j_h\,\Div\Big(\sum_i
y_i\bfphi^i_h\Big)\,dx\quad\text{for every }j,
\end{align*}
whence, $\sum_i y_i\bfphi^i_h\in \mathscr P^{k, h}_{0, \Div}(\Omega,
\rn)$. Thus, owing to (\ref{cor:4.1.0}),
\begin{align*}
\langle\bfb,\bfy\rangle= \int_\Omega \bfH_h:\nabla
\Big(\sum_iy_i\bfphi^i_h\Big)\,dx=0
\end{align*}
for every $\bfy \in \Ker(\mathcal A^T)$. Hence,  $\bfb\in
\big(\Ker\mathcal A^T\big)^\perp$.
\par\noindent Next,  we show that
$\mathcal A$ is injective, whence the uniqueness of $\pi_h$ follows.
To verify the injectivity of $\mathcal A$, assume that $\bfz$ is
such that $\mathcal A \bfz=0$. Hence, on setting $p_h =\sum_j z_j\,
p_h^j$, we have that $\int_\Omega p_h\,\Div\bfphi^i_h\,dx=0$ for
every $i$.
Let us apply  Theorem \ref{thm:4.2}  with any pair of Young
functions $A$ and $B$ fulfilling \eqref{1.1} and \eqref{1.2}, for
instance $A(t)=B(t)=t^2$ for $t \geq 0$ (this means that we are in
fact applying \eqref{infsupp} with $p=2$). We then obtain that
\begin{align*}
\|p_h\|_{L^{2}(\Omega)}\leq C\sup_{\bfvarphi_h\in \mathscr P^{k,
h}_0(\Omega, \rn)}\frac{\int_\Omega
p_h\,\Div\bfphi_h\,dx}{\|\nabla \bfvarphi_h\|_{L^2(\Omega,
\rmat)}}=0,
\end{align*}
%
whence  $\bfz=0$.
 \\
(ii) By Theorem \ref{thm:4.2}
 and (\ref{B.5})
\begin{align*}
\|\pi_h\|_{L^B(\Omega)}&\leq C \sup_{\bfvarphi_h\in \mathscr P^{k,
h}_0(\Omega, \rn)}\frac{\int_\Omega
\pi_h\,\Div\bfphi_h\,dx}{\|\nabla
\bfvarphi_h\|_{L^{\widetilde{A}}(\Omega, \rmat)}} =C
\sup_{\bfvarphi_h\in \mathscr P^{k, h}_0(\Omega,
\rn)}\frac{\int_\Omega
\bfH_h:\nabla\bfphi_h\,dx}{\|\nabla
\bfvarphi_h\|_{L^{\widetilde{A}}(\Omega, \rmat)}}\\
&\leq C \sup_{\bfvarphi\in W^{1,\widetilde{A}}_0(\Omega,
\rn)}\frac{\int_\Omega \bfH_h:\nabla\bfphi\,dx}{\|\nabla
\bfvarphi\|_{L^{\widetilde{A}}(\Omega, \rmat)}} \leq C
\sup_{\bfPsi\in L^{\widetilde{A}}(\Omega, \rmat)}\frac{\int_\Omega
\bfH_h:\bfPsi\,dx}{\|\bfPsi\|_{L^{\widetilde{A}}(\Omega, \rmat)}}\\
&\leq C'  \|\bfH_h\|_{L^A(\Omega, \rmat)}
\end{align*}
for some constants $C=C(\Omega, c)$ and $C'=C'(\Omega , c)$. This
proves inequality \eqref{oct1}.
\\ (iii) The triangle inequality and the embedding
$L^A(\Omega)\rightarrow L^B(\Omega)$ ensure that, for every
$\mu_h\in\mathscr P^{m, h}_\bot(\Omega)$,
\begin{align*}
\|\pi_h-\pi\|_{L^B(\Omega)}&\leq \|\pi_h-\mu_h\|_{L^B(\Omega)}+\|\mu_h-\pi\|_{L^B(\Omega)}\\
&\leq \|\pi_h-\mu_h\|_{L^B(\Omega)}+C\|\mu_h-\pi\|_{L^A(\Omega)},
\end{align*}
for some constant $C=C(c)$. By Theorem \ref{thm:4.2},
\begin{align*}
\|\pi_h-\mu_h\|_{L^B(\Omega)}&\leq C \sup_{\bfvarphi_h\in \mathscr
P^{k, h}_0(\Omega, \rn)}\frac{\int_\Omega
(\pi_h-\mu_h)\,\Div\bfphi_h\,dx}{\|\nabla
\bfvarphi_h\|_{L^{\widetilde{A}}(\Omega, \rmat)}}
\end{align*}
for some constant $C(\Omega , c)$, and for every $\mu_h\in\mathscr
P^{m, h}_\bot(\Omega)$.
 Moreover, by
\eqref{march1}, \eqref{march2}, and an approximation argument for
functions in $\mathscr P^{k, h}_{0}(\Omega, \rn)$ via functions in
$C^\infty _0(\Omega, \rn)$, we have that
\begin{align*}
\int_\Omega (\pi_h-\mu_h)\,\Div\bfphi_h\,dx&=\int_\Omega
(\pi_h-\pi)\,\Div\bfphi_h\,dx+\int_\Omega
(\pi-\mu_h)\,\Div\bfphi_h\,dx\\
&=\int_\Omega (\bfH_h-\bfH):\nabla\bfphi_h\,dx+\int_\Omega
(\pi-\mu_h)\,\Div\bfphi_h\,dx,
\end{align*}
for every $\mu_h\in\mathscr P^{m, h}_\bot(\Omega)$ and $\bfvarphi_h
\in \mathscr P^{k, h}_{0}(\Omega, \rn)$.
 Hence, via H\"older's inequality in Orlicz spaces,
\begin{align*}
\int_\Omega (\pi_h-\mu_h)\,\Div\bfphi_h\,dx&\leq C
\bigg(\|\bfH_h-\bfH\|_{L^A(\Omega,
\rmat)}+\|\mu_h-\pi\|_{L^A(\Omega)}\bigg)\|\nabla\bfphi_h\|_{L^{\widetilde{A}}(\Omega,
\rmat)}
\end{align*}
for some constant $C=C(n)$. Altogether, we conclude that
\begin{align*}
\|\pi_h-\pi\|_{L^B(\Omega)}\leq C\bigg(\|\bfH_h-\bfH\|_{L^A(\Omega,
\rmat)}+\inf_{\mu_h\in\mathscr P^{m,
h}_\bot(\Omega)}\|\mu_h-\pi\|_{L^A(\Omega)}\bigg)
\end{align*}
for some constant $C$ independent of $h$, namely \eqref{oct2}.
\qed

\bigskip
\par\noindent
{\bf Acknowledgements}. We wish to thank M. Bulicek for some useful
discussions. \par\noindent This research was partly supported  by
Leopoldina (German National Academy of Science),  by the Research
Project  ``Geometric properties of partial differential equations
and related topics" 2008 of the Italian MIUR (Ministry of Education,
University and Research), and by GNAMPA of the Italian INdAM
(Istituto Nazionale di Alta Matematica).

\end{document}